\documentclass[12pt]{article}

\usepackage{fontenc,fontaxes}
\usepackage{amsthm}
\usepackage{amssymb}
\usepackage{lmodern}
\usepackage{authblk}
\usepackage[authoryear,round]{natbib}
\usepackage{comment}
\usepackage{mathtools}
\usepackage{tikz}
\usepackage[bf]{caption}
\usepackage{mathabx}
\usepackage{breakcites}
\usepackage{setspace}
\usepackage{microtype}
\usepackage[greek,english]{babel}

\usepackage[breaklinks=true,bookmarksnumbered,bookmarksopen,colorlinks,linkcolor={black},urlcolor={black},citecolor={black}]{hyperref}
\newtheorem{thm}{Theorem}[section]
\newtheorem{cor}{Corollary}[section]

\newtheorem*{AFT}{Analytic Fredholm Theorem}

\theoremstyle{definition}
\newtheorem{remark}{Remark}[section]
\newtheorem{assumption}{Assumption}[section]
\newtheorem{definition}{Definition}[section]
\newtheorem{example}{Example}[section]

\mathtoolsset{mathic=true}

\DeclareMathOperator{\ran}{ran}

\DeclareMathOperator{\coker}{coker}
\DeclareMathOperator{\cl}{cl}
\DeclareMathOperator{\spn}{sp}

\newcommand{\id}{\mathrm{id}}
\newcommand{\zero}{0}
\numberwithin{equation}{section}

\renewcommand\footnotemark{}

\begin{document}
	
	\newcommand{\gpi}{\textrm{\greektext p}}	
	
	\title{Representation of I(1) and I(2) autoregressive Hilbertian processes}
	
\author[1]{Brendan K.\ Beare}
\author[2]{Won-Ki Seo}\thanks{We thank Massimo Franchi, Peter Phillips and seminar participants at the Einaudi Institute for Economics and Finance, UC San Diego, the Universit\'{e} libre de Bruxelles, and the 2019 NBER-NSF Time Series Conference in Hong Kong for helpful discussions. Beare also thanks Phil Roberts for his feedback and encouragement. An earlier version of this article titled ``Representation of I(1) autoregressive Hilbertian processes'' was posted on the arXiv.org preprint repository in January 2017.}
\affil[1]{School of Economics, University of Sydney}
\affil[2]{Department of Economics, Queen's University}
	\maketitle
\begin{center}
	Accepted for publication in \textit{Econometric Theory}.
\end{center}	
	\bigskip
	\begin{abstract}	
		We develop versions of the Granger-Johansen representation theorems for I(1) and I(2) vector autoregressive processes that apply to processes taking values in an arbitrary complex separable Hilbert space. This more general setting is of central relevance for statistical applications involving functional time series. An I(1) or I(2) solution to an autoregressive law of motion is obtained when the inverse of the autoregressive operator pencil has a pole of first or second order at one. We obtain a range of necessary and sufficient conditions for such a pole to be of first or second order. Cointegrating and attractor subspaces are characterized in terms of the behavior of the autoregressive operator pencil in a neighborhood of one.		
		
		
		
		
		
	\end{abstract}
	
	\pagebreak
	\onehalfspacing
	
	\section{Introduction}\label{intro}

	Results on the existence and representation of integrated solutions to vector autoregressive laws of motion are among the most important and subtle contributions of econometricians to time series analysis, yet also among the most widely misunderstood. The best known such result is the so-called Granger representation theorem, which first appeared in an unpublished UC San Diego working paper of \citet{G83}. In this paper, Granger, having recently introduced the concept of cointegration \citep{G81} sought to connect statistical models of time series based on linear process representations to regression-based models involving equilibrium correction terms, which had appeared in work by \citet{S64} and \citet{DHSY78}. The main result of \citet{G83} first emerged in published form in \citet{G86} without proof, but more prominently in the widely cited \emph{Econometrica} article by \citet{EG87}, where it is labeled the ``Granger representation theorem'', with the exclusion of the first author presumably due to the paper having resulted from the merger of previous independent contributions.
	
	The proof of the Granger representation theorem in \citet{EG87} is incorrect. Moreover, the error can be traced back to the original working paper of \citet{G83}. A counterexample to Lemma A1 of \citet{EG87}, which is also Theorem 1 of \citet{G83}, may be found buried in a footnote of \citet{J09}. Johansen was familiar with Granger's work on representation theory at an early stage, visiting UC San Diego and authoring a closely related Johns Hopkins working paper in 1985 that was eventually published as \citet{J88}. At around the same time the doctoral thesis of \citet{Y87} at UC San Diego established the connection to Smith-McMillan forms. \citet{J91} provided what appears to be the first correct statement and proof of a modified version of the Granger representation theorem, which we will call the Granger-Johansen representation theorem. This contribution did not merely correct a technical error of Granger; it reoriented attention toward a central issue: when does a given vector autoregressive law of motion admit an I(1) solution? The answer to this question is given by the Johansen I(1) condition, which is a necessary and sufficient condition on the autoregressive polynomial and its first derivative at one for a vector autoregressive law of motion to admit an I(1) solution. \citet{J92} subsequently obtained analogous results for I(2) processes, and \citet{L98} provided some results for higher-order integrated processes.
	
	An early contribution of \citet{S91} contained a striking observation on the Johansen I(1) condition: it corresponds to a necessary and sufficient condition for the inverse of a matrix pencil to have a simple pole at a given point in the complex plane. Various authors later exploited this insight, including \citet{FZ02,FZ09,FZ11}, \citet{F07}, \citet{J09} and \citet{FP11,FP16}. A nice aspect of the connection to analytic function theory is that it extends naturally to the development of \(\mathrm{I}(d)\) representation theory with integer-valued \(d\geq2\): just as the Johansen I(1) condition can be reformulated as a necessary and sufficient condition for a simple pole, analogous \(\mathrm{I}(d)\) conditions  can be reformulated as necessary and sufficient conditions for poles of order \(d\). \citet{FP19} have recently taken precisely this approach to develop a general \(\mathrm{I}(d)\) representation theory. The introduction to their paper contains a detailed discussion of the history of research on the Granger-Johansen representation theorem.
	
	Parallel to the development of representation theorems for cointegrated systems in the 1980s and early 1990s was the development of asymptotic distribution theory for the statistical estimation of such systems, obtained by applying central limit theory on function spaces and associated results. This research was led by Phillips and his students at Yale; see, in particular, \citet{PD86}, \citet{P86,P88,P91}, \citet{PPk88}, \citet{PkP88,PkP89} and \citet{PH90}, among many other contributions. \citet{J91} used limit theorems developed in this body of work, applicable to integrated systems exhibiting general serial dependence, to derive the asymptotic distribution of the maximum likelihood estimator for I(1) Gaussian vector autoregressive systems. Complementary work by \citet{CW87,CW88} and \citet{KP91} on asymptotic theory for regression with integrated processes has also been influential.
	
	In this paper we provide representation theorems for I(1) and I(2) vector autoregressive processes taking values in an arbitrary complex separable Hilbert space. This more general setting is of central relevance for statistical applications involving functional time series \citep{HK12}, and was first studied by \citet{CKP16} in the case of I(1) probability density-valued time series; see also \citet{B17} and \citet{SB19}. Our results here build on those we obtained in an earlier paper with J.\ Seo \citep{BSS17} establishing a representation theorem for the I(1) case. They differ from our earlier results in their explicit use of analytic function theory. In Theorems \ref{mainthm} and \ref{mainthm2}, our I(1) and I(2) representation theorems, we directly assume that the inverse of the autoregressive operator pencil has a pole of first or second order at one. We elaborate upon the meaning of these pole conditions in Theorems \ref{polethm} and \ref{polethm2}, which provide necessary and sufficient conditions to have a pole of first or second order. These results also provide explicit formulas for the coefficients in the principal part of the corresponding Laurent series.
	
	Our paper supersedes an earlier manuscript posted on the arXiv.org preprint repository in January 2017 \citep{BS17} that dealt only with the I(1) case. During its preparation several working papers have emerged that deliver related results. In particular, \citet{FP18} study I($d$) solutions to autoregressive laws of motion in complex separable Hilbert space, for integer-valued $d\geq1$. Their necessary and sufficient condition for an I($d)$ solution involves an orthogonal direct sum decomposition of the Hilbert space into $d$ closed subspaces. This contrasts with the direct sum conditions given by \citet{BSS17} for the I(1) case, and here for the I(1) and I(2) cases, which involve nonorthogonal direct sums. We also provide a range of alternative formulations of our necessary and sufficient conditions, some of which may be easier to verify than others. Also relevant is recent work by \citet{HP17}, who established an alternative I($d$) condition for first-order autoregressive Hilbertian processes: the restriction of the autoregressive operator to the image of the Riesz projection associated with its unit eigenvalue differs from the identity by an operator nilpotent of degree $d$. Finally, \citet{CHP16} have developed I(1) representation theory for autoregressive Hilbertian processes under the assumption that the impact operator in the error correction representation is compact. Under this condition the dimension of the cointegrating space must be finite, which contrasts with the setting of this paper and the others cited in this paragraph, where the codimension of the cointegrating space must be finite. Finite codimensionality of the cointegrating space implies that the I(1) stochastic trend in the Beveridge-Nelson representation of our cointegrated process is confined to a finite dimensional linear subspace. It is a consequence of a compactness condition we impose on the autoregressive operators. \citet{FP18} have observed that finite codimensionality of the cointegrating space holds more generally if the autoregressive operator pencil has an eigenvalue of finite type at one.
	
	We structure the remainder of the paper as follows. Section \ref{ssetup} sets the scene with notation and essential mathematics. Our results on I(1) and I(2) processes are contained in Sections \ref{sI1} and \ref{sI2} respectively. We provide a brief discussion of directions for future research in Section \ref{sconc}. Appendix \ref{aftapx} contains background material on the spectral properties of operator-valued functions, including a statement of the analytic Fredholm theorem, which is a key input to our results. The proofs of our results are collected in Appendix \ref{appx}.

	\section{Preliminaries}\label{ssetup}
	
	The setting for our analysis is a separable complex Hilbert space \(\mathcal H\) with inner product \(\langle\cdot{,}\cdot\rangle\) and norm \(\Vert\cdot\Vert\). If \(\mathcal H'
	\) is another such space, we let \(\mathcal L_{\mathcal H,\mathcal H'}\) denote the Banach space of continuous linear operators from \(\mathcal H\) to \(\mathcal H'\) equipped with the operator norm. We are mostly concerned with the case \(\mathcal H=\mathcal H'\), and write \(\mathcal L_{\mathcal H}\) in place of \(\mathcal L_{\mathcal H,\mathcal H}\). To each operator \(A\in\mathcal{L}_{\mathcal H}\) we associate two linear subspaces of \(\mathcal H\): the kernel and range of \(A\), given by
	\begin{equation*}
	\ker A=\{x\in \mathcal H:A(x)=0\},\quad\ran A=\{A(x):x\in \mathcal H\}.
	\end{equation*}
	We let \(\mathrm{I}\in\mathcal{L}_{\mathcal H}\) denote the identity map on \(\mathcal H\).
	
	A central concern of our analysis will be the decomposability of \(\mathcal H\) into sums of certain linear subspaces of \(\mathcal H\). Given linear subspaces $V$ and $W$ of $\mathcal H$, we write $V+W$ for the linear subspace of all $x\in\mathcal H$ such that \(x=v+w\) for some $v\in V$ and $w\in W$. When $V$ and $W$ are linear subspaces of $\mathcal{H}$ with \(V\cap W=\{0\}\), we may instead write $V\oplus W$ for their sum, and call it a direct sum. When we write \(\mathcal H=V\oplus W\), we are asserting that \(\mathcal H=V+W\) and that \(V\cap W=\{0\}\). In this case, any \(x\in\mathcal H\) may be uniquely decomposed as \(x=v+w\) with \(v\in V\) and \(w\in W\).
	
	Orthogonal complements and projections play a key role in our analysis. Given a linear subspace \(V\) of \(\mathcal H\), we define its orthogonal complement by
	\begin{equation*}
	V^\perp=\{x\in\mathcal H:\langle x,v\rangle=0\text{ for all }v\in V\}.
	\end{equation*}
	The orthogonal complement to a linear subspace of \(\mathcal H\) is always a closed linear subspace of \(\mathcal H\). Given a closed linear subspace \(V\) of \(\mathcal H\), it is always the case that \(\mathcal H=V\oplus V^\perp\). Thus any \(x\in\mathcal H\) may be uniquely decomposed as \(x=v+v'\) with \(v\in V\) and \(v'\in V^\perp\). We denote by \(\mathrm{P}_V\in\mathcal{L}_{\mathcal H}\) the orthogonal projection on \(V\), which maps a point \(x=v+v'\) to \(v\).
	
	Our main results concern the representation of time series taking values in \(\mathcal H\), but only the most basic understanding of probability on \(\mathcal H\) is required. As in \citet{BSS17}, we let \(L_{\mathcal H}^2\) denote the Banach space of random elements \(Z\) of \(\mathcal H\) (identifying random elements that are equal with probability one) that satisfy \(E\Vert Z\Vert^2<\infty\) and \(EZ=0\), equipped with the norm \(\Vert Z\Vert_{L^2_{\mathcal H}}=(E\Vert Z\Vert^2)^{1/2}\). Refer to that paper for the definition of \(EZ\) and of the covariance operator of an element of \(L_{\mathcal H}^2\). For further details, the monograph of \citet{B00} provides a comprehensive treatment of linear processes taking values in a real Hilbert or Banach space. A complex Hilbert space setting was studied more recently by \citet{CH17}.

	\section{I(1) autoregressive Hilbertian processes}\label{sI1}
	
	In this section we state our results for I(1) autoregressive processes. Let \(p \in \mathbb{N}\), and consider the following AR(\(p\)) law of motion in \(\mathcal{H} \):
	\begin{equation}\label{arp}
	X_t=\sum_{j=1}^p\Phi_j(X_{t-j})+\varepsilon_t.
	\end{equation}
	Here, the \(X_t\)'s and \(\varepsilon_t\)'s are random elements of $\mathcal{H}$, and the \(\Phi_j\)'s are continuous linear operators from \(\mathcal{H}\) to \(\mathcal{H}\). We say that the AR(\(p\)) law of motion \eqref{arp} is engendered by the map \(\Phi : \mathbb{C} \mapsto \mathcal L_\mathcal{H} \) given by
	\begin{align}\label{arpencil}
	\Phi(z) = \mathrm{I} - \sum_{j=1}^p z^j \Phi_j.
	\end{align}
	We will refer to an operator-valued polynomial function of a complex variable as an operator pencil; note that some authors reserve this term for linear polynomials. We impose the following conditions on the objects just introduced.
	\begin{assumption}\label{assume1}
		(i) \(\varepsilon = (\varepsilon_t, t \in \mathbb{Z})\) is an iid sequence in \( L^2_\mathcal{H}\) with positive definite covariance operator \( \Sigma \in \mathcal L_\mathcal{H}\). (ii) \( \Phi_1, \ldots, \Phi_p\) are compact operators in \(\mathcal L_\mathcal{H}\) such that \( \Phi : \mathbb{C} \mapsto \mathcal L_\mathcal{H}\) is noninvertible at \(z=1 \) and invertible at every other \(z\) in the closed unit disk.
	\end{assumption}
	\begin{remark}
		The innovations \(\varepsilon_t\) are referred to as strong white noise due to their being centered (i.e.\ zero expected value) and iid. We impose these conditions for simplicity, but the results to be developed remain valid if the iid condition is replaced with the weaker requirement that the cross-covariance operators for the \(\varepsilon_t\)'s are all zero, as in \citet{FP18}. In the latter case the \(\varepsilon_t\)'s are merely said to be white noise. More generally, one might consider allowing the \(\varepsilon_t\)'s to be a general I(0) process as in \citet{J92} and \citet{CP09,CP12}, or even a subexponential process as in \citet{A18}, but we do not pursue this route here.
	\end{remark}
	\begin{remark}
		The results to be developed remain valid if the assumption that \(\Phi_1,\ldots,\Phi_p\) are compact is replaced with the weaker but less easily interpretable requirement that \(\Phi(z)\) has an eigenvalue of finite type at \(z=1\). See \citet{FP18} for details.
	\end{remark}
	
	 The approach we will take to developing representation theory for I(1) and I(2) autoregressive processes in \(\mathcal H\) essentially boils down to studying the behavior of \(\Phi(z)^{-1}\) near \(z=1\). We achieve this by applying the analytic Fredholm theorem, a complete statement of which is provided in Appendix \ref{aftapx}. To apply this result we need \(\Phi(z)\) to be analytic in \(z\) (in fact, it is polynomial in \(z\), hence analytic), and Fredholm operator-valued (meaning that it has finite dimensional kernel and cokernel for all $z$, which is guaranteed by our compactness condition on the \(\Phi_j\)'s). The analytic Fredholm theorem implies that, for all \(z\) in a punctured neighborhood of one, we have
	\begin{equation}\label{aftlaurentintro}
	\Phi(z)^{-1}=\sum_{k=-d}^\infty(z-1)^k\Upsilon_k,
	\end{equation}
	where \(d\in\mathbb N\) and \(\Upsilon_{-d},\Upsilon_{-d+1},\ldots\) is a sequence in \(\mathcal L_\mathcal{H}\). The series in \eqref{aftlaurentintro} is called the Laurent series of \(\Phi(z)^{-1}\) around \(z=1\), and converges in \(\mathcal L_{\mathcal H}\). If we assume without loss of generality that \(\Upsilon_{-d}\neq\zero\), then \(\Phi(z)^{-1}\) is said to have a pole of order \(d\) at \(z=1\). The operator \(\Upsilon_{-1}\) is called the residue of \(\Phi(z)^{-1}\) at \(z=1\). A pole of order one is said to be simple. We call the sum of the leading terms indexed by \(k=-d,\ldots,-1\) the principal part of the Laurent series, and we call the truncated series excluding these leading terms the analytic part of the Laurent series.
	
	An important implication of the analytic Fredholm theorem is that the leading Laurent coefficients \(\Upsilon_{-d},\ldots,\Upsilon_{-1}\) in \eqref{aftlaurentintro} are all of finite rank. In the representation theory to be developed, this has the effect of ensuring that we always obtain a cointegrating space with finite codimension (i.e., with finite dimensional orthogonal complement).
	
	It will be convenient to introduce some additional notation. We will write \(\Pi_0\) and \(\Pi_1\) for the values taken by \(\Phi(z)\) and its first derivative at \(z=1\):
	\begin{equation*}
	\Pi_0=\Phi(1),\quad\Pi_1=\Phi^{(1)}(1).
	\end{equation*}
	We also define the linear spaces
	\begin{equation*}
	\alpha_1=\ran\Pi_0,\quad\beta_1=(\ker\Pi_0)^\perp.
	\end{equation*}
	Our compactness condition on the \(\Phi_j\)'s ensures that \(\alpha_1\) and \(\beta_1\) are closed linear spaces with equal and finite codimension.
	
	We have yet to give a formal definition of the I(\(d\)) property. For our purposes, it is sufficient to define the I(0) property for standard linear processes. We will need to consider standard linear processes in \(\mathcal H\) and in \(\mathbb C\), with innovations in \(\mathcal H\), so give the following definition of a standard linear process in an arbitrary separable complex Hilbert space \(\mathcal H'\).
	\begin{definition}
		A sequence \((W_t,t\geq t_0)\) in \(L^2_{\mathcal H'}\) is called a standard linear process in \(\mathcal H'\) if there is another separable complex Hilbert space \(\mathcal H\) such that we may write
		\begin{equation}\label{linproc}
		W_t=\sum_{k=0}^\infty A_k(\varepsilon_{t-k}),\quad t\geq t_0,
		\end{equation}
		where \((A_k,k\geq0)\) is a norm-summable sequence in \(\mathcal{L}_{\mathcal H, \mathcal H'}\), and \((\varepsilon_t,t\in\mathbb Z)\) is an iid sequence in \(L^2_{\mathcal H}\) with nonzero covariance operator \(\Sigma\in\mathcal L_{\mathcal H}\).
	\end{definition}
\begin{remark}\label{iplprem}
	If we were to require that the two Hilbert spaces \(\mathcal H\) and \(\mathcal H'\) are the same, and that \(A_0=\mathrm{I}\), then our definition of a standard linear process in \(\mathcal H\) would be the same as that of \citet[p.\ 183]{B00}. Our more general definition is needed because if \((W_t,t\geq t_0)\) is a standard linear process in \(\mathcal H\) with innovations in \(\mathcal H\), then for any \(x\in\mathcal H\) we may write
	\begin{equation*}
	\langle x,W_t\rangle=\sum_{k=0}^\infty \langle x,A_k(\varepsilon_{t-k})\rangle=\sum_{k=0}^\infty \hat{A}^x_k(\varepsilon_{t-k}),\quad t\geq t_0,
	\end{equation*}
	where the \(\hat{A}^x_k\)'s are given by
	\begin{equation*}
	\hat{A}^x_k(y)=\langle x,A_k(y)\rangle,\quad y\in\mathcal H,
	\end{equation*}
	and form a norm-summable sequence in \(\mathcal{L}_{\mathcal H, \mathbb C}\). The sequence of inner products \((\langle x,W_t\rangle,t\geq t_0)\) is thus a standard linear process in \(\mathbb C\) with innovations in \(\mathcal H\).
\end{remark}
We may now define the I(\(d\)) property for sequences in \(L^2_{\mathcal H'}\), with \(\mathcal H'\) an arbitrary separable complex Hilbert space.
	\begin{definition}\label{I0def}
		We say that a sequence \((W_t,t\geq t_0)\) in \(L^2_{\mathcal H'}\) is I(0) if it is a standard linear process in \(\mathcal H'\) admitting a representation \eqref{linproc} in which \(\Sigma\) is positive definite and the \(A_k\)'s satisfy \(\sum_{k=0}^\infty A_k\neq0\) and \(\sum_{k=0}^\infty k\Vert A_k\Vert_{\mathcal L_{\mathcal H,\mathcal H'}}<\infty\).
	\end{definition}
	\begin{definition} We say that a sequence \((W_t,t\geq t_0)\) in \(L^2_{\mathcal H'}\) is I(\(d\)) for \(d\in\mathbb N\) if its \(d\)th difference is an I(0) standard linear process in \(\mathcal H'\).
		\end{definition}
		\begin{remark}
			The summability condition on the norms of the coefficients \(A_k\) in Definition \ref{I0def} is called 1-summability. It was used by \citet{PS92} to facilitate a version of the Beveridge-Nelson decomposition for a time series whose difference is I(0). In the results to be developed, all processes claimed to be I(0) in fact have coefficients decaying exponentially in norm, so 1-summability is easily satisfied.
		\end{remark}
	
	Our first result provides an I(1) representation for autoregressive Hilbertian processes for which \(\Phi(z)^{-1}\) has a simple pole at \(z=1\). We will discuss the simple pole condition in more detail later in this section.

	\begin{thm}\label{mainthm}
		Suppose that Assumption \ref{assume1} is satisfied, and that \(\Phi(z)^{-1}\) has a simple pole at \(z=1\). Let \(\Upsilon_{-1}\) denote the residue of \(\Phi(z)^{-1}\) at \(z=1\), let \(\tilde{\Psi}(z) \) denote the analytic part of the Laurent series of \( \Phi(z)^{-1} \) around \(z=1\), and set \( \tilde{\Psi}_k = \tilde{\Psi}^{(k)}(0)/k! \). A sequence \( (X_t, t \geq -p+1) \) in \(L^2_\mathcal{H}\) satisfying the law of motion \eqref{arp} allows the following representation: for some \(Z_0\in L^2_\mathcal{H} \) and all \(t\geq1\) we have
		\begin{equation}\label{deltax}
		X_t = Z_0 -\Upsilon_{-1}\left(\sum_{s=1}^t \varepsilon_{s}\right) + \nu_t .
		\end{equation}
		Here,  \( (\nu_t,t\geq1)\) is a stationary sequence of random elements of \(\mathcal{H}\) defined by the \(L^2_\mathcal{H}\)-convergent series \(\nu_t=\sum_{k=0}^\infty \tilde{\Psi}_k (\varepsilon_{t-k}) \). Moreover,
		\begin{itemize}
			\item[\((1)\)] The range of \(\Upsilon_{-1}\) is equal to \(\beta_1^\perp\) and has positive and finite dimension;
			\item[\((2)\)] If \(Z_0\) belongs to \(\beta_1^\perp\), then for nonzero \(x\in \mathcal{H}\) the sequence of inner products \((\langle x,X_t\rangle,t\geq1)\) is \(\mathrm{I}(0)\) if \(x\in\beta_1\), and is \(\mathrm{I}(1)\) otherwise.
		\end{itemize}
	\end{thm}
		\begin{remark}
			\hyperref[mainthm]{Theorem \ref*{mainthm}} is similar to Theorem 4.1 of \citet{BSS17}, but makes the connection to the analytic behavior of \(\Phi(z)^{-1}\) explicit. The latter result is more generally applicable in one respect: compactness of the autoregressive operator is not assumed when \(p=1\). The approach taken here relies on the analytic Fredholm theorem and therefore requires \(\Phi(z)\) to be Fredholm, which may not be the case if the autoregressive operators are not compact.
		\end{remark}
		\begin{remark}\label{finite}
			The residue \(\Upsilon_{-1}\) appearing in \hyperref[mainthm]{Theorem \ref*{mainthm}} has finite rank by the analytic Fredholm theorem. The attractor space, which is the subspace of \(\mathcal{H}\) in which the I(1) stochastic trend in the Beveridge-Nelson representation \eqref{deltax} takes values, thus has finite dimension. We are therefore outside the framework considered by \citet{CHP16}, in which the cointegrating space has finite dimension and the attractor space has finite codimension.
		\end{remark}
		
		When can we expect the simple pole condition in Theorem \ref{mainthm} to be satisfied? Our next result provides equivalent reformulations of this condition that may be easier to check in practice, and a formula for the residue \(\Upsilon_{-1}\) in terms of \(\alpha_1\), \(\beta_1\) and \(\Pi_1\). (Recall we defined \(\Pi_1=\Phi^{(1)}(1)\).)
		
	\begin{thm}\label{polethm}
			Suppose that Assumption \ref{assume1}(ii) holds. The following four conditions are equivalent.
		\begin{itemize}
			\item[\((1)\)] \(\Phi(z)^{-1}\) has a simple pole at \(z=1\).
			\item[\((2)\)] The operator \(\Lambda_1:\beta_1^\perp\to\alpha_1^\perp\) obtained by restricting \(\mathrm{P}_{\alpha_1^\perp}\Pi_1\) to \(\beta_1^\perp\) is bijective.
			\item[\((3)\)] \(\mathcal{H}=\alpha_1\oplus\Pi_1\beta_1^\perp\).
			\item[\((4)\)] \(\mathcal{H}=\alpha_1+\Pi_1\beta_1^\perp\).
		\end{itemize}
	If \(\Phi(z)^{-1}\) has a simple pole at \(z=1\), then its residue at \(z=1\) is \(\Upsilon_{-1}=\Lambda_1^{-1}\mathrm{P}_{\alpha_1^\perp}\).
	\end{thm}

		\begin{remark}
			The closest results we have found to \hyperref[polethm]{Theorem \ref*{polethm}} in prior literature are those of \citet{S68} and \citet{H71}, who worked in a more general Banach space setting. \citet{S68} established sufficient conditions for a simple pole, and \citet{H71} established the equivalence of conditions (1) and (3).
		\end{remark}
		\begin{remark}\label{remJI1}
			\citet[Prop.\ 4.1]{BSS17} showed that condition (3) of Theorem \ref{polethm} is equivalent to the I(1) condition given by \citet[Thm.\ 4.1]{J91} in the finite dimensional case \(\mathcal H=\mathbb C^n\). In this setting we may let \(r<n\) be the rank of the \(n\times n\) complex matrix \(\Pi_0\), let \(\alpha\) and \(\beta\) be full-rank \(n\times r\) complex matrices such that \(\Pi_0=\alpha\beta^\prime\), and let \(\alpha_\perp\) and \(\beta_\perp\) be full-rank \(n\times(n-r)\) complex matrices such that \(\alpha^\prime\alpha_\perp=0\) and \(\beta^\prime\beta_\perp=0\). The Johansen I(1) condition is satisfied when the \((n-r)\times(n-r)\) complex matrix \(\alpha_\perp^\prime\Pi_1\beta_\perp\) is invertible. Examples 4.1--4.3 of \citet{BSS17} illustrate the reformulation of the Johansen I(1) condition as a direct sum decomposition.
		\end{remark}
		\begin{remark}
			The direct sum appearing in condition (3) of Theorem \ref{polethm} is not in general an orthogonal direct sum. \citet{FP18} showed that, when \(\Phi(z)\) is noninvertible at \(z=1\), condition (3) is equivalent to the following orthogonal direct sum decomposition of \(\mathcal H\):
			\[\mathcal H=\beta_1\oplus\left(\ker\Lambda_1\mathrm{P}_{\beta_1^\perp}\right)^\perp. \]
			In their notation, this is \(\mathcal H=\tau_0\oplus\tau_1\).
		\end{remark}
		\begin{remark}
			Our assumption that the operators \(\Phi_1,\ldots,\Phi_p\) are compact implies that \(\Pi_0\) is Fredholm of index zero. If \(\Pi_0\) were Fredholm but not of index zero then it would be impossible to satisfy condition (2). This is because bijectivity of \(\Lambda_1\) requires its domain and codomain to have the same dimension. However, our proof that condition (4) implies condition (1) does not use the index-zero property, and remains valid if the compactness condition on \(\Phi_1,\ldots,\Phi_p\) is weakened to require only that \(\Phi(z)\) is Fredholm operator-valued.
		\end{remark}
		In the special case where \(p=1\), conditions (3) and (4) of \hyperref[polethm]{Theorem \ref*{polethm}} take on a particularly simple form, and another related equivalent condition becomes available. Moreover, the direct sum decomposition asserted by condition (3) serves to define an oblique projection that is the negative of the residue of our simple pole. The following corollary to \hyperref[polethm]{Theorem \ref*{polethm}} provides details.

	\begin{cor}\label{polethmp1}
		Suppose that Assumption \ref{assume1}(ii) holds and that \(p=1\). The following four conditions are equivalent.
		\begin{itemize}
			\item[\((1)\)] \(\Phi(z)^{-1}\) has a simple pole at \(z=1\).
			\item[\((2)\)] \(\mathcal{H}=\alpha_1\oplus\beta_1^\perp\).
			\item[\((3)\)] \(\mathcal{H}=\alpha_1+\beta_1^\perp\).
			\item[\((4)\)] \(\{0\}=\alpha_1\cap\beta_1^\perp\).
		\end{itemize}
	If \(\Phi(z)^{-1}\) has a simple pole at \(z=1\), then its residue at \(z=1\) is the negative of the projection on \(\beta_1^\perp\) along \(\alpha_1\).
	\end{cor}

	\begin{remark}\label{rieszrem}
	The oblique projection appearing in \hyperref[polethmp1]{Corollary \ref*{polethmp1}} is in fact the Riesz projection for the unit eigenvalue of \(\Phi_1\). Said Riesz projection is defined (\citealp[p.\ 9]{GGK90}; \citealp[pp.\ 11--12]{M12}) by the contour integral
	\begin{equation}\label{riesz}
	P=\frac{1}{2\gpi\mathrm{i}}\oint_\Gamma(z\mathrm{I}-\Phi_1)^{-1}\mathrm{d}z,
	\end{equation}
	where \(\Gamma\) is a positively oriented smooth Jordan curve around one separating it from zero and from any other eigenvalues of \(\Phi_1\), and where the integral of an \(\mathcal L_\mathcal{H}\)-valued function should be understood in the sense of Bochner. Let \(\gamma:[0,1]\to\mathbb C\) be a smooth parametrization of \(\Gamma\), and rewrite \eqref{riesz} as
	\begin{equation}\label{riesz2}
	P=\frac{1}{2\gpi\mathrm{i}}\int_0^1(\gamma(t)\mathrm{I}-\Phi_1)^{-1}\gamma'(t)\mathrm{d}t.
	\end{equation}
	The image of \(\Gamma\) under the reciprocal transform \(z\mapsto z^{-1}\), which we denote \(\Gamma'\), is a positively oriented smooth Jordan curve around one separating it from any other poles of \(\Phi(z)^{-1}\) and from zero. It admits the parametrization \(t\mapsto1/\gamma(t)\eqqcolon\delta(t)\). A little calculus shows that \(\gamma'(t)=-\delta'(t)/\delta(t)^2\), and so from \eqref{riesz2} we have
	\begin{equation}\label{riesz3}
	P=\frac{-1}{2\gpi\mathrm{i}}\int_0^1\delta(t)^{-1}(\mathrm{I}-\delta(t)\Phi_1)^{-1}\delta'(t)\mathrm{d}t
	=\frac{-1}{2\gpi\mathrm{i}}\oint_{\Gamma'} z^{-1}\Phi(z)^{-1}\mathrm{d}z.
	\end{equation}
	The residue theorem therefore tells us that \(P\) is the negative of the residue of \(z^{-1}\Phi(z)^{-1}\) at \(z=1\), implying that the residue of \(\Phi(z)^{-1}\) at \(z=1\) is \(-P\). It now follows from \hyperref[polethmp1]{Corollary \ref*{polethmp1}} that when the direct sum decomposition \(\mathcal{H}=\alpha_1\oplus\beta_1^\perp\) is satisfied, the Riesz projection for the unit eigenvalue of \(\Phi_1\) is the projection on \(\beta_1^\perp\) along \(\alpha_1\). This is a nonorthogonal projection except when \(\alpha_1=\beta_1\), which would occur if, for instance, \(\Phi_1\) is a normal operator.
\end{remark}
\begin{remark}
	It is apparent from our discussion in Remark \ref{rieszrem} that, under any of the equivalent I(1) conditions in Corollary \ref{polethmp1}, the restriction of \(\Pi_0\) to the range of the Riesz projection \(P\) is zero. This is precisely the condition given by \citet{HP17} for an AR(1) process in \(\mathcal H\) to be I(1).
\end{remark}
We close this section with several examples of the use of Corollary \ref{polethmp1} for verifying that \(\Phi(z)^{-1}\) has a simple pole at \(z=1\).

\begin{example} \label{ex1}
	Suppose that \(p=1\) and that \(\Phi_1\) is self-adjoint. Then \(\Pi_0\) is also self-adjoint, implying that \(\alpha_1=\beta_1\). The direct sum decomposition \(\mathcal{H} = \alpha_1 \oplus \beta_1^\perp\) appearing in Corollary \ref{polethmp1} is therefore satisfied, and is in this case an orthogonal direct sum decomposition.    
\end{example}

\begin{example} \label{ex2}
	Let \((e_j , j \in \mathbb{N})\) be an orthonormal basis of \(\mathcal{H}\). Suppose that \(p=1\), and that \(\Phi_1\) is given by
	\begin{equation*}
	\Phi_1(x) = \langle x, e_1\rangle (e_1+e_2) +  \sum_{j=2}^\infty \lambda_j \langle x, e_j\rangle e_j,\quad x\in \mathcal{H},
	\end{equation*}
	with $\lambda_j\in(0,1)$ for \(j\geq2\) and $\lambda_j \rightarrow 0$ as $j \rightarrow \infty$. For any \(x \in \mathcal{H}\) with representation \(x= \sum_{j=1}^\infty c_j e_j\), \(c_j=\langle x,e_j\rangle\), we have
	\begin{align}
	(\mathrm{I}-\Phi_1)(x) = (c_2(1-\lambda_2)-c_1)e_2 + \sum_{j=3}^\infty  c_j(1-\lambda_j)e_j. \label{exameq1}
	\end{align}
	Since \(\lambda_j\neq1\) for all \(j\geq3\), it is clear that \(e_j\notin\ker (\mathrm{I}-\Phi_1)\) for all \(j\geq3\). Moreover,
	\begin{align*}
	(\mathrm{I}-\Phi_1)(c_1 e_1 +c_2 e_2) = (c_2 (1-\lambda_2) - c_1) e_2.
	\end{align*}
	It follows that
	\begin{align}\label{ex2ker}
	\beta_1^\perp=\ker (\mathrm{I}-\Phi_1)=\{c_1e_1+c_2e_2:c_1=c_2(1-\lambda_2)\}.
	\end{align}
	Moreover, it may be deduced that $\alpha_1=\cl\spn\{e_j:j\geq2\}$, the closed linear span of \(\{e_j:j\geq2\}\), as follows. Any \(x\in\cl\spn\{e_j:j\geq2\}\) may be written as \(x = \sum_{j=2}^\infty d_je_j\) for some square-summable sequence \( (d_j, j\geq2) \). We can always find another square-summable sequence \((c_j, j \in \mathbb{N})\) such that 
	\begin{align}
	d_2=c_2(1-\lambda_2) - c_1 \quad \text{and} \quad  d_j = c_j(1-\lambda_j),\quad j\geq3.
	\end{align} 
	Then
	\begin{align*}
	(\mathrm{I}-\Phi_1)\left(\sum_{j=1}^\infty c_je_j\right)=(c_2(1-\lambda_2)-c_1)e_2 + \sum_{j=3}^\infty  c_j(1-\lambda_j)e_j=\sum_{j=2}^\infty d_je_j=x,
	\end{align*}
	which shows that $x\in\ran(\mathrm{I}-\Phi_1)=\alpha_1$. Thus \(\cl\spn\{e_j:j\geq2\}\subseteq\alpha_1 \). In addition, it is easily deduced that \(\alpha_1 \subseteq  \cl\spn\{e_j:j\geq2\} \) using \eqref{exameq1}. Therefore, \(\alpha_1 =\cl\spn\{e_j:j\geq2\}\). From \eqref{ex2ker} we see that the only element of \(\beta_1^\perp\) belonging to \(\cl\spn\{e_j:j\geq2\}\) is zero. Thus the condition \(\{0\}=\alpha_1\cap\beta_1^\perp\) appearing in \hyperref[polethmp1]{Corollary \ref*{polethmp1}} is satisfied.
\end{example}

\begin{example} \label{ex3}
	Suppose that in Example \ref{ex2} we instead defined \(\Phi_1\in\mathcal L_\mathcal{H}\) by   
	\begin{equation*}
	\Phi_1(x) = \langle x, e_1\rangle (e_1+e_2+e_3) + \langle x, e_2\rangle e_2  + \langle x, e_3\rangle e_3 + \sum_{j=4}^\infty \lambda_j \langle x, e_j\rangle e_j,\quad x\in \mathcal{H},
	\end{equation*}
	with $\lambda_j\in(0,1)$ for $j\geq4$ and $\lambda_j \rightarrow 0$ as $j \rightarrow \infty$. For any \(x \in \mathcal{H}\) with representation \(x= \sum_{j=1}^\infty c_j e_j\), \(c_j=\langle x,e_j\rangle\), we now have
	\begin{align*}
	(\mathrm{I}-\Phi_1)(x) =& -c_1e_2 -c_1e_3+ \sum_{j=4}^\infty c_j(1-\lambda_j)e_j.
	\end{align*}
	Since \(\lambda_j\neq1\) for all \(j\geq4\), it is clear that \(e_j \notin \ker (\mathrm{I}-\Phi_1) \) for all \(j\geq4\). Moreover,
	\begin{align}\label{exop6}
	(\mathrm{I}-\Phi_1)(c_1 e_1 +c_2 e_2 +c_3 e_3) =  -c_1 e_2 - c_1 e_3. 
	\end{align}
	It follows that $\beta_1^\perp=\ker (\mathrm{I}-\Phi_1)=\spn \{e_2, e_3\}$. Further, arguments similar to those in Example \ref{ex2} can be used to show that 
	\begin{align*}
	\alpha_1 = \cl\spn\{e_2+e_3, e_4, e_5, \ldots\}.
	\end{align*}
	It follows that \(\alpha_1\cap\beta_1^\perp=\spn\{e_2+e_3\}\). Thus the condition \(\{0\}=\alpha_1\cap\beta_1^\perp\) appearing in \hyperref[polethmp1]{Corollary \ref*{polethmp1}} is violated.	
\end{example}

\section{I(2) autoregressive Hilbertian processes}\label{sI2}

In this section we state our results for I(2) autoregressive processes. We continue to assume that the autoregressive coefficients and innovations of such a process satisfy Assumption \ref{assume1}. It will be convenient to introduce some additional notation. We define the linear spaces
\begin{equation}\label{alpha2beta2}
\alpha_2=\alpha_1 + \Pi_1\beta_1^\perp,\quad
\beta_2=\beta_1+\Pi_1^\ast\alpha_1^\perp.
\end{equation}
The role played by \(\alpha_2\) and \(\beta_2\) in our I(2) results will be analogous to the role played by \(\alpha_1\) and \(\beta_1\) in our I(1) results. The notation \(\Pi_1^\ast\) refers to the adjoint operator to \(\Pi_1\).

\begin{remark}\label{kercokerrem}
	The linear spaces \(\alpha_2^\perp\) and \(\beta_2^\perp\) have equal and finite dimension and are, respectively, the cokernel and kernel of the operator \(\Lambda_1\) appearing in Theorem \ref{polethm}. To see why, observe that \(\alpha_2\) satisfies
	\begin{equation*}
	\alpha_2^\perp=(\alpha_1+\mathrm{P}_{\alpha_1^\perp}\Pi_1\beta_1^\perp)^\perp=(\alpha_1+\ran\Lambda_1)^\perp=\ker\Lambda_1^\ast,
	\end{equation*}
	and similarly \(\beta_2\) satisfies
	\begin{equation*}
	\beta_2^\perp=(\beta_1+\mathrm{P}_{\beta_1^\perp}\Pi_1^\ast\alpha_1^\perp)^\perp=(\beta_1+\ran\Lambda_1^\ast)^\perp=\ker\Lambda_1.
	\end{equation*}
	The operator \(\Lambda_1\) is invertible under any of the equivalent I(1) conditions given in Theorem \ref{polethm}, in which case we must have \(\alpha_2=\beta_2=\mathcal H\). In this section we are interested in the case where our I(1) conditions fail. This occurs when \(\alpha_2^\perp\) and \(\beta_2^\perp\) have positive dimension.
\end{remark}

\begin{remark}\label{ab2p1rem}
Simpler expressions for \(\alpha_2\) and \(\beta_2\) become available when \(p=1\). In this case we have \(\Pi_1=\Pi_0-\mathrm{I}\), from which it follows easily that \(\Pi_1\beta_1^\perp=\beta_1^\perp\) and \(\Pi_1^\ast\alpha_1^\perp=\alpha_1^\perp\). We may therefore write
\begin{equation*}
\alpha_2=\alpha_1 + \beta_1^\perp,\quad
\beta_2=\beta_1+\alpha_1^\perp.
\end{equation*}
\end{remark}

Our first result in this section provides an I(2) analogue to Theorem \ref{mainthm} in Section \ref{sI1}. It establishes an I(2) representation for autoregressive Hilbertian processes for which \(\Phi(z)^{-1}\) has a pole of second order at \(z=1\). We will discuss the pole condition in more detail later in this section.

\begin{thm}\label{mainthm2}
	Suppose that Assumption \ref{assume1} is satisfied, and that \(\Phi(z)^{-1}\) has a pole of second order at \(z=1\). Let \(\Upsilon_{-2}\) and \(\Upsilon_{-1}\) denote the coefficients in the principal part of the Laurent series of \( \Phi(z)^{-1} \) around \(z=1\), let \(\tilde{\Psi}(z) \) denote the analytic part of the Laurent series of \( \Phi(z)^{-1} \) around \(z=1\), and set \( \tilde{\Psi}_k = \tilde{\Psi}^{(k)}(0)/k! \). A sequence \( (X_t, t \geq -p+1) \) in \(L^2_\mathcal{H}\) satisfying the law of motion \eqref{arp} allows the following representation: for some \(Z_0,Z_1\in L^2_\mathcal{H} \) and all \(t\geq1\) we have
	\begin{equation}\label{deltax2}
	X_t = Z_0+tZ_1 + \Upsilon_{-2}\left(\sum_{s=1}^t\sum_{r=1}^s \varepsilon_{r}\right)-\Upsilon_{-1}\left(\sum_{s=1}^t\varepsilon_{s}\right) + \nu_t .
	\end{equation}
	Here,  \( (\nu_t,t\geq1)\) is a stationary sequence of random elements of \(\mathcal{H}\) defined by the \(L^2_\mathcal{H}\)-convergent series \(\nu_t = \sum_{k=0}^\infty \tilde{\Psi}_k (\varepsilon_{t-k}) \). Moreover,
	\begin{itemize}
		\item[\((1)\)] The range of \(\Upsilon_{-2}\) is equal to \(\beta_2^\perp\) and has positive and finite dimension;
		\item[\((2)\)] If \(Z_0\) and \(Z_1\) belong to \(\beta_2^\perp\), then for nonzero \(x\in \mathcal{H}\) the sequence of inner products \((\langle x,X_t\rangle,t\geq1)\) is \(\mathrm{I}(0)\) or \(\mathrm{I}(1)\) if \(x\in\beta_2\), and is \(\mathrm{I}(2)\) otherwise.
	\end{itemize}
\end{thm}

	In view of claim (2) in Theorem \ref{mainthm2} we may refer to the linear subspace \(\beta_2\) as the cointegrating space. For \(x\in\beta_2\), the sequence of inner products \((\langle x,\Delta X_t\rangle,t\geq1)\) may be \(\mathrm{I}(0)\) or \(\mathrm{I}(1)\) (ignoring the effect of the deterministic components \(Z_0\) and \(Z_1\)). Polynomial cointegration may also occur. This will be discussed more fully in Remark \ref{FPI2} below.

Our next result provides an I(2) analogue to Theorem \ref{polethm} in Section \ref{sI1}. It establishes necessary and sufficient conditions for \(\Phi(z)^{-1}\) to have a pole of second order at \(z=1\), as assumed in Theorem \ref{mainthm2}. It also establishes formulas for the Laurent coefficients \(\Upsilon_{-2}\) and \(\Upsilon_{-1}\) appearing in Theorem \ref{mainthm2}. To state our result it will be convenient to introduce some additional notation. Similar to \(\Pi_0\) and \(\Pi_1\), we define
\begin{equation*}
\Pi_2=\frac{1}{2}\Phi^{(2)}(1)-\Pi_1\Pi_0^+\Pi_1,\quad
\Pi_3=\frac{1}{6}\Phi^{(3)}(1)-\Pi_1\Pi_0^+\Pi_1\Pi_0^+\Pi_1.
\end{equation*}
Here and elsewhere, a plus superscript is used to denote the Moore-Penrose inverse of a continuous linear operator between Hilbert spaces with closed range \citep[see e.g.][ch.\ 9]{BG03}.

\begin{thm}\label{polethm2}
	Suppose that Assumption \ref{assume1}(ii) holds, and that \(\Phi(z)^{-1}\) does not have a simple pole at \(z=1\). Then the following four conditions are equivalent.
	\begin{itemize}
		\item[\((1)\)] \(\Phi(z)^{-1}\) has a pole of second order at \(z=1\).
		\item[\((2)\)] The operator \(\Lambda_2: \beta_2^\perp \to \alpha_2^\perp\) obtained by restricting \(\mathrm{P}_{\alpha_2^\perp}\Pi_2\) to \(\beta_2^\perp\) is bijective.
		\item[\((3)\)] \(\mathcal{H}=\alpha_2 \oplus \Pi_2 \beta_2^\perp\).
		\item[\((4)\)] \(\mathcal{H}=\alpha_2+ \Pi_2 \beta_2^\perp\).
	\end{itemize}
	If \(\Phi(z)^{-1}\) has a pole of second order at \(z=1\), then the coefficients of \((z-1)^{-2}\) and \((z-1)^{-1}\) in the Laurent series of \( \Phi(z)^{-1} \) around \(z=1\) are given by
	\begin{equation}\label{resdef2}
		\Upsilon_{-2}=\Lambda_2^{-1}\mathrm{P}_{\alpha_2^\perp}
	\end{equation}
	and
	\begin{align}\label{resdef22}
	\Upsilon_{-1}&=\Lambda_1^+  \mathrm{P}_{\alpha_1^\perp}-\left[\Lambda_1^+  \mathrm{P}_{\alpha_1^\perp}\Pi_2+\Pi_0^+ \Pi_1\right]\Upsilon_{-2}-\Upsilon_{-2}\left[\Pi_2\Lambda_1^+  \mathrm{P}_{\alpha_1^\perp}+\Pi_1\Pi_0^+\right]\notag\\
	&\quad+\Upsilon_{-2}\left[\Pi_1\Pi_0^+\Pi_2+\Pi_2\Pi_0^+\Pi_1+\Pi_2\Lambda_1^+\mathrm{P}_{\alpha_1^\perp}\Pi_2-\Pi_3\right]\Upsilon_{-2},
	\end{align}
	respectively.
\end{thm}
	\begin{remark}
	The four equivalent conditions in Theorem \ref{polethm2} are equivalent to the I(2) condition given by \citet[Thm.\ 3]{J92} in the finite dimensional case \(\mathcal H=\mathbb C^n\). Continuing with the notation of Remark \ref{remJI1}, suppose that the Johansen I(1) condition fails, and let \(\varphi\) and \(\eta\) be full-rank \((n-r)\times s\) (\(s<n-r\)) complex matrices such that \((\alpha_\perp^\prime\alpha_\perp)^{-1}\alpha_\perp^\prime\Pi_1\beta_\perp(\beta_\perp^\prime\beta_\perp)^{-1}=\varphi\eta^\prime\). Let \(\tilde{\alpha}_\perp\) and \(\tilde{\beta}_\perp\) be full-rank \(n\times(n-r-s)\) complex matrices whose columns are orthogonal to those of \(\tilde{\alpha}\coloneqq(\alpha,\alpha_\perp\varphi)\) and \(\tilde{\beta}\coloneqq(\beta,\beta_\perp\eta)\) respectively. Note that the column spaces of \(\tilde{\alpha}\) and \(\tilde{\beta}\) are, respectively, the linear spaces \(\alpha_2\) and \(\beta_2\) defined in \eqref{alpha2beta2} above. The Johansen I(2) condition is satisfied when the \((n-r-s)\times(n-r-s)\) complex matrix \(\tilde{\alpha}_\perp^\prime\Pi_2\tilde{\beta}_\perp\) is invertible. Lemma 4.1 of \citet{BSS17} implies that \(\tilde{\alpha}_\perp^\prime\Pi_2\tilde{\beta}_\perp\) is invertible if and only if \(\mathcal H\) is the direct sum of the null space of \(\tilde{\alpha}_\perp^\prime\) and the column space of \(\Pi_2\tilde{\beta}_\perp\). The former space is \(\alpha_2\) and the latter space is \(\Pi_2\beta_2^\perp\), so the Johansen I(2) condition is equivalent to condition (3) in Theorem \ref{polethm2}.
    \end{remark}
	\begin{remark}
		Formulas \eqref{resdef2} and \eqref{resdef22} in Theorem \ref{polethm2} correspond to formulas (12) and (13) of \citet[Thm.\ 5]{J09}, given for the finite dimensional case \(\mathcal H=\mathbb C^n\). The objects \(C_1\), \(C_2\), \(\bar{\beta}_1\bar{\alpha}_1^\prime\), \(\bar{\beta}\bar{\alpha}'\), \(\dot{\Pi}\), \(\theta\) and \((1/6)\dddot{\Pi}-\dot{\Pi}\bar{\beta}\bar{\alpha}'\dot{\Pi}\bar{\beta}\bar{\alpha}'\dot{\Pi}\) in Johansen's notation correspond respectively to \(-\Upsilon_{-1}\), \(\Upsilon_{-2}\), \(\Lambda_1^+\mathrm{P}_{\alpha_1^\perp}\), \(-\Pi_0^+\), \(\Pi_1\), \(\Pi_2\) and \(\Pi_3\) in our notation.
	\end{remark}
	\begin{remark}\label{FPI2}
		The direct sum appearing in condition (3) of Theorem \ref{polethm2} is not in general an orthogonal direct sum. Extending results of \citet{J92} from the finite dimensional case \(\mathcal H=\mathbb C^n\) to a more general Hilbert space setting, \citet{FP18} showed that, when \(\Phi(z)\) is noninvertible at \(z=1\) and the I(1) condition fails, an equivalent necessary and sufficient condition for a pole of second order is the following tripartite orthogonal direct sum decomposition of \(\mathcal H\):
		\begin{equation*}
		\mathcal H=\beta_1\oplus\left(\ker\Lambda_1\mathrm{P}_{\beta_1^\perp}\right)^\perp\oplus\left(\ker\Lambda_2\mathrm{P}_{\beta_2^\perp}\right)^\perp.
		\end{equation*}
		In their notation, this is \(\mathcal H=\tau_0\oplus\tau_1\oplus\tau_2\). Moreover, we have the following bipartite orthogonal direct sum decomposition of the cointegrating space \(\beta_2\):
		\begin{equation*}
		\beta_2=\beta_1\oplus\left(\ker\Lambda_1\mathrm{P}_{\beta_1^\perp}\right)^\perp,
		\end{equation*}
		or \(\beta_2=\tau_0\oplus\tau_1\). This decomposition is informative because the two subspaces decomposing the cointegrating space correspond to different kinds of cointegrating behavior. Suppose that \(Z_0\) and \(Z_1\) belong to \(\beta_2^\perp\), so that we may ignore the effect of deterministic components. For \(x\in\beta_2\) not belonging to \(\beta_1\), the sequence \((\langle x,X_t\rangle,t\geq1)\) is \(\mathrm{I}(1)\). For nonzero \(x\in\beta_1\), the sequence \((\langle x,X_t\rangle,t\geq1)\) may be \(\mathrm{I}(0)\) or \(\mathrm{I}(1)\), but the sequence \((\langle x,X_t\rangle-\langle x,\Pi_0^+\Pi_1\Delta X_t\rangle,t\geq1)\) is always \(\mathrm{I}(0)\). The case where the inner product sequences \((\langle x,X_t\rangle,t\geq1)\) and \((\langle x,\Pi_0^+\Pi_1\Delta X_t\rangle,t\geq1)\) are \(\mathrm{I}(1)\) and cointegrated is called polynomial cointegration or multicointegration \citep{Y87,GL89,GL90,EJ99,PK19}. 
	\end{remark}
	\begin{remark}\label{I0space}
		From our discussion of the \(\mathrm{I}(2)\) case in Remark \ref{FPI2}, it is apparent that (ignoring deterministic components) the sequence \((\langle x,X_t\rangle,t\geq1)\) is \(\mathrm{I}(0)\) for precisely those nonzero \(x\in\beta_1\) for which \((\langle x,\Pi_0^+\Pi_1\Delta X_t\rangle,t\geq1)\) is stationary. This is the collection of all nonzero \(x\in\beta_1\) such that \((\Pi_0^+\Pi_1)^\ast(x)\) belongs to \(\beta_2\), the cointegrating space. Therefore, \((\langle x,X_t\rangle,t\geq1)\) is \(\mathrm{I}(0)\) if and only if \(x\) is a nonzero element of the linear space
		\begin{equation*}
		\zeta=\beta_1\cap(\Pi_0^+\Pi_1\beta_2^\perp)^\perp.
		\end{equation*}
	\end{remark}	
	\begin{remark}\label{pole2p1rem}
		When \(p=1\), the linear space \(\Pi_2\beta_2^\perp\) appearing in conditions (3) and (4) of Theorem \ref{polethm2} is equal to \((\mathrm{I}-\Pi_0^+)(\alpha_1\cap\beta_1^\perp)\). To see why, observe that when \(p=1\) we have
		\begin{equation*}
		\Pi_2=-\Pi_1\Pi_0^+\Pi_1=-(\mathrm{I}-\Pi_0)\Pi_0^+(\mathrm{I}-\Pi_0)=-\Pi_0^++\Pi_0^+\Pi_0+\Pi_0\Pi_0^+-\Pi_0\Pi_0^+\Pi_0.
		\end{equation*}
		From Remark \ref{ab2p1rem} we know that \(\beta_2^\perp=\alpha_1\cap\beta_1^\perp\) when \(p=1\). Therefore, since \(\Pi_0\) is zero on \(\beta_1^\perp\), we have
		\begin{equation*}
		\Pi_2\beta_2^\perp=(\Pi_0\Pi_0^+-\Pi_0^+)(\alpha_1\cap\beta_1^\perp).
		\end{equation*}
		From the basic properties of Moore-Penrose inverses, we know that \(\Pi_0\Pi_0^+\) is orthogonal projection on \(\alpha_1\), the range of \(\Pi_0\). The restriction of this projection to \(\alpha_1\cap\beta_1^\perp\) coincides with the identity, so our claim is established.
	\end{remark}	
	\begin{remark}\label{I0spacep1}
		When \(p=1\), the linear space \(\zeta\) appearing in Remark \ref{I0space} is equal to \(\beta_1\cap(\Pi_0^+(\alpha_1\cap\beta_1^\perp))^\perp\). To see why, observe that when \(p=1\) we have \(\beta_2^\perp=\alpha_1\cap\beta_1^\perp\) (shown in Remark \ref{ab2p1rem}) and \(\Pi_0^+\Pi_1=\Pi_0^+\Pi_0-\Pi_0^+\). From the basic properties of Moore-Penrose inverses, we know that \(\Pi_0^+\Pi_0\) is orthogonal projection on \(\beta_1\), the corange of \(\Pi_0\). Thus \(\Pi_0^+\Pi_0\) is zero on \(\alpha_1\cap\beta_1^\perp\), and we have \(\Pi_0^+\Pi_1\beta_2^\perp=\Pi_0^+(\alpha_1\cap\beta_1^\perp)\), establishing our claim.
	\end{remark}

			\begin{remark}\label{geneigenrem}
				As discussed in Remark \ref{rieszrem}, when \(p=1\), the negative of the residue \(\Upsilon_{-1}\) is a projection, and is called the Riesz projection for the unit eigenvalue of \(\Phi_1\). The space on which \(-\Upsilon_{-1}\) projects is called the generalized eigenspace for the unit eigenvalue, and its dimension is called the algebraic multiplicity of the unit eigenvalue. Contained within the generalized eigenspace is the usual eigenspace \(\beta_1^\perp\), whose dimension is called the geometric multiplicity of the unit eigenvalue \citep[p.\ 26]{GGK90}. The generalized eigenspace is in fact the orthogonal complement to \(\zeta\), the subspace of the cointegrating space \(\beta_2\) yielding \(\mathrm{I}(0)\) inner products, as defined in Remark \ref{I0space}. To see why, observe that \(\zeta^\perp\) is the subspace occupied by the \(\mathrm{I}(1)\) and \(\mathrm{I}(2)\) trends in the representation \eqref{deltax2}, which is the sum of the ranges of \(\Upsilon_{-2}\) and \(\Upsilon_{-1}\). But from formula \eqref{resdef2} in Theorem \ref{polethm2}, we know that the range of \(\Upsilon_{-2}\) is \(\beta_2^\perp\), which is contained in the usual eigenspace \(\beta_1^\perp\), and therefore contained in the generalized eigenspace, which is the range of \(\Upsilon_{-1}\). Therefore \(\zeta^\perp\) is the generalized eigenspace. In view of Remark \ref{I0spacep1}, we have
				\begin{equation}\label{geneigenspacep1}
				\zeta^\perp=\beta_1^\perp+\Pi_0^+(\alpha_1\cap\beta_1^\perp).
				\end{equation}
				From Corollary \ref{polethmp1} we know that for \(p=1\) the \(\mathrm{I}(1)\) condition fails precisely when \(\alpha_1\cap\beta_1^\perp\neq\{0\}\). Since the Moore-Penrose inverse \(\Pi_0^+\) defines a bijection from \(\alpha_1\) to \(\beta_1\), and \(\beta_1^\perp\) is finite dimensional, we deduce that when \(p=1\) and the \(\mathrm{I}(1)\) condition fails we must have \(\dim\beta_1^\perp<\dim\zeta^\perp.\) Thus we see that when the \(\mathrm{I}(1)\) condition fails, the algebraic multiplicity of the unit eigenvalue exceeds its geometric multiplicity. This contrasts with the situation when the \(\mathrm{I}(1)\) condition is satisfied, where, as is apparent from our discussion in Remark \ref{rieszrem}, the algebraic and geometric multiplicities of the unit eigenvalue are equal. The fact that the equality of geometric and algebraic multiplicities of the unit eigenvalue implies an \(\mathrm{I}(1)\) representation was observed by \citet[Cor.\ 4.3]{J96} in the finite dimensional case \(\mathcal H=\mathbb C^n\).
			\end{remark}
	\begin{remark}
		\citet{HP17} provide the following condition for an AR(1) process in \(\mathcal H\) to be \(\mathrm{I}(2)\): the restriction of \(\Pi_0\) to the generalized eigenspace \(\zeta^\perp\) must be nilpotent of degree two. We can verify this condition using the expression we obtained for \(\zeta^\perp\) in \eqref{geneigenspacep1}. Based on this expression, any element of \(\zeta^\perp\) may be written as \(x+\Pi_0^+(y)\), where \(x\in\beta_1^\perp\) and \(y\in\alpha_1\cap\beta_1^\perp\). Observe that
		\begin{equation*}
		\Pi_0(x+\Pi_0^+(y))=\Pi_0\Pi_0^+(y)=y\quad\text{and}\quad\Pi_0^2(x+\Pi_0^+(y))=\Pi_0(y)=0,
		\end{equation*}
		since \(\Pi_0\Pi_0^+\) is orthogonal projection on \(\alpha_1\). Thus the restriction of \(\Pi_0\) to \(\zeta^\perp\) is zero if and only if \(\alpha_1\cap\beta_1^\perp=\{0\}\), which is one of our necessary and sufficient conditions for a simple pole from Corollary \ref{polethmp1}. Assuming that this condition fails, and that our pole is of second order, we find that the restriction of \(\Pi_0\) to \(\zeta^\perp\) is nilpotent of degree two.
	\end{remark}

\begin{figure}[t!]
	\begin{center}
		\begin{tikzpicture}[scale=3]
		\draw[black, thick, <->] (0,-1) -- (0,1);
		\draw[black, thick, <->] (-1,0) -- (1,0);
		\draw[black] (0,.05) -- (.05,.05) -- (.05,0);
		\draw[black,thick,domain=-1:1, <->]  plot(\x, \x);
		\node[below] at (1,0) {\footnotesize \(x_1\)};
		\node[right] at (0,1) {\footnotesize \(x_2\)};
		\node[below,rotate=90] at (-.2,0.5) {\scriptsize \(\alpha_1=\beta_1^\perp=\alpha_2=\beta_2^\perp\)};
		\node[below] at (.5,0.15) {\scriptsize \(\beta_1=\beta_2\)};
		\node[above, rotate=45] at (.58,.55) {\scriptsize \(\Pi_2\beta_2^\perp\)};
		\node[below] at (0,-1.2) {(a) Example \ref{jex2}};
		\end{tikzpicture}
		\begin{tikzpicture}[scale=3]
		\draw[black, thick, <->] (0,-1) -- (0,1);
		\draw[black, thick, <->] (-1,0) -- (1,0);
		\draw[black,thick,domain=-1:1, <->]  plot(\x, \x);
		\draw[black,thick,domain=-1:1, <->]  plot(\x, -\x);
		\draw[black,thick,domain=-1:1, <->]  plot(\x, -3*\x/5);
		\draw[black] (0,.05) -- (.05,.05) -- (.05,0);
		\node[below] at (1,0) {\footnotesize \(x_1\)};
		\node[left] at (0,1) {\footnotesize \(x_2\)};
		\node[below, rotate=-45] at (.655,-.6) {\scriptsize \(\alpha_1=\beta_1=\Pi_1\beta_1^\perp=\alpha_2=\beta_2\)};
		\node[above, rotate=45] at (.58,.55) {\scriptsize \(\beta_1^\perp=\beta_2^\perp\)};
		\node[above, rotate={atan(-3/5)}] at (.675,-.45) {\scriptsize \(\Pi_2\beta_2^\perp\)};
		\node[below] at (0,-1.2) {(b) Example \ref{jex3}};
		\end{tikzpicture}
	\end{center}
	\caption{Visual aid to Examples \ref{jex2}--\ref{jex3}.}
	\label{jexfig}
\end{figure}

To illustrate the preceding results on \(\mathrm{I}(2)\) representations we discuss three examples. The first two examples correspond to Examples 4.2--4.3 of \citet{J96} and \citet{BSS17}, with \(\gamma=2\) in the latter example so that the solution is I(2). These examples are finite (two) dimensional.

\begin{example}\label{jex2}
	Consider the first-order autoregressive law of motion \(X_t=\Phi_1(X_{t-1})+\varepsilon_t\) in \(\mathcal H=\mathbb C^2\) with autoregressive coefficient matrix
	\begin{equation*}
	\Phi_1=\left[\begin{array}{cc}1&0\\1&1\end{array}\right].
	\end{equation*}
	The only eigenvalue of \(\Phi_1\) is one. With a little algebra we find that
	\begin{equation*}
	\alpha_1=\beta_1^\perp=\alpha_2=\beta_2^\perp=\mathrm{sp}\left[\begin{array}{c}0\\1\end{array}\right],\quad\Pi_2\beta_2^\perp=\mathrm{sp}\left[\begin{array}{c}1\\1\end{array}\right].
	\end{equation*}
	We depict these subspaces of \(\mathbb C^2\) in Figure \ref{jexfig}(a), displaying them as subspaces of \(\mathbb R^2\) since \(\Phi_1\) has real elements. (We do the same for Example \ref{jex3} in Figure \ref{jexfig}(b).) It is apparent that the (version for \(p=1\) of the) I(1) condition \(\mathbb C^2=\alpha_1\oplus\beta_1^\perp\) is not satisfied, and that the I(2) condition \(\mathbb C^2=\alpha_2\oplus\Pi_2\beta_2^\perp\) is satisfied, so our autoregressive law of motion generates an I(2) process. It is easily verified that \(\zeta=\{0\}\), so in this case it follows from Theorem \ref{mainthm2} that (ignoring deterministic components) the inner product sequence \((\langle x,X_t\rangle,t\geq1)\) is I(1) for nonzero \(x\in\beta_2\) and I(2) for nonzero \(x\notin\beta_2\).
\end{example}

\begin{example}\label{jex3}
	Consider the second-order autoregressive law of motion \(X_t=\Phi_1(X_{t-1})+\Phi_2(X_{t-2})+\varepsilon_t\) in \(\mathcal H=\mathbb C^2\) with autoregressive coefficient matrices
	\begin{equation*}
	\Phi_1=\left[\begin{array}{cc}\frac{5}{4}&\frac{7}{4}\\-\frac{1}{4}&\frac{5}{4}\end{array}\right],\quad\Phi_2=\left[\begin{array}{cc}0&-2\\0&0\end{array}\right].
	\end{equation*}
	The points of noninvertibility of \(\Phi(z)\) are \(z=1\) and \(z=2\). It is straightforward to establish that
	\begin{equation*}
	\alpha_1=\beta_1=\Pi_1\beta_1^\perp=\alpha_2=\beta_2=\mathrm{sp}\left[\begin{array}{c}1\\-1\end{array}\right],\quad\Pi_2\beta_2^\perp=\mathrm{sp}\left[\begin{array}{c}5\\-3\end{array}\right].
	\end{equation*}
	We depict these subspaces in Figure \ref{jexfig}(b). It is apparent that the I(1) condition \(\mathbb C^2=\alpha_1\oplus\Pi_1\beta_1^\perp\) is not satisfied, and that the I(2) condition \(\mathbb C^2=\alpha_2\oplus\Pi_2\beta_2^\perp\) is satisfied, so our autoregressive law of motion generates an I(2) process. As in Example \ref{jex2}, we have \(\zeta=\{0\}\), so that (ignoring deterministic components) \((\langle x,X_t\rangle,t\geq1)\) is I(1) for nonzero \(x\in\beta_2\) and I(2) for nonzero \(x\notin\beta_2\).
\end{example}

Our final example is infinite dimensional, and builds on Example \ref{ex3} above.
\begin{example} \label{ex4}
	Consider the setting of Example \ref{ex3}. We can use Theorem \ref{polethm2} to determine whether we have a pole of second order. The Moore-Penrose inverse \(\Pi_0^+\) satisfies \(\Pi_0^+\Pi_0=\mathrm{P}_{\beta_1}\). Applying \(\Pi_0^+\) to both sides of the equality \((\mathrm{I}-\Phi_1)(-e_1)=e_2+e_3\) reveals that \(\mathrm{P}_{\beta_1}(-e_1)=\Pi_0^+(e_2+e_3)\), which simplifies to \(\Pi_0^+(e_2+e_3)=-e_1\) since \(\beta_1=\spn\{e_2,e_3\}^\perp\). In view of Remark \ref{pole2p1rem}, it follows that
	\begin{equation*}
	\Pi_2\beta_2^\perp=(\mathrm{I}-\Pi_0^+)(\alpha_1\cap\beta_1^\perp)=(\mathrm{I}-\Pi_0^+)\spn\{e_2+e_3\}=\spn\{e_1+e_2+e_3\}.
	\end{equation*}
	In view of Remark \ref{ab2p1rem}, we have
	\begin{equation*}
	\alpha_2=\alpha_1+\beta_1^\perp=\cl\spn\{e_2+e_3,e_4,e_5,\ldots\}+\spn\{e_2,e_3\}=\cl\spn\{e_2,e_3,e_4,\ldots\}.
	\end{equation*}
	Thus we see that \(\mathcal{H}\) is the sum of the linear subspaces \(\alpha_2\) and \(\Pi_2\beta_2^\perp\), and deduce from Theorem \ref{polethm2} that \(\Phi(z)^{-1}\) has a pole of second order at \(z=1\). The associated cointegrating space is
	\begin{equation*}
	\beta_2=\beta_1+\alpha_1^\perp=\spn\{e_2,e_3\}^\perp+\cl\spn\{e_2+e_3,e_4,e_5,\ldots\}^\perp=\cl\spn\{e_1,e_2-e_3,e_4,\ldots\}.
	\end{equation*}
	The \(\mathrm{I}(2)\) stochastic trend takes values in the orthogonal complement to this space, which is
	\begin{equation*}
	\beta_2^\perp=\spn\{e_2+e_3\}.
	\end{equation*}
	In view of Remark \ref{geneigenrem}, the \(\mathrm{I}(1)\) stochastic trend takes values in the larger space
	\begin{equation*}
	\zeta^\perp=\beta_1^\perp+\Pi_0^+(\alpha_1\cap\beta_1^\perp)=\spn\{e_2,e_3\}+\spn\{e_1\}=\spn\{e_1,e_2,e_3\}.
	\end{equation*}
	It follows that (ignoring deterministic components) for nonzero \(x\in\mathcal H\) the inner product sequence \((\langle x,X_t\rangle,t\geq1)\) is
	\begin{align*}
	\text{I(0) for }&x\in\cl\spn\{e_4,e_5,\ldots\},\\
	\text{I(1) for }&x\in\cl\spn\{e_1,e_2-e_3,e_4,\ldots\}\setminus\cl\spn\{e_4,e_5,\ldots\},\\
	\text{I(2) for }&x\notin\cl\spn\{e_1,e_2-e_3,e_4,\ldots\}.
	\end{align*}

	
\end{example}

\section{Concluding remarks}\label{sconc}
	
	The results established in this paper extend long-known results on I(1) and I(2) representations for autoregressive processes in \(\mathbb C^n\) to a more general Hilbert space setting. It may be desirable to extend our results further to a Banach space setting. In a Banach space setting we no longer have the luxury of using inner products to formulate a suitable notion of cointegration. It is instead natural to define cointegration in terms of the order of integration of continuous linear functionals of an integrated process, and the cointegrating space as a subspace of the topological dual. An unpublished manuscript of the second author \citep{S18} provides some results in this direction, and also investigates relaxing the compactness condition we have imposed here on autoregressive operators. It is interesting to drop compactness because in this case the attractor space may be infinite dimensional; the difficulty is that the analytic Fredholm theorem becomes unavailable.
	
	Beyond representation theory, research on the development of statistical procedures for analyzing cointegrated functional time series is a priority. A new manuscript by \citet{NSS19} develops a procedure for estimating the dimension of an attractor space based on sequential variance ratio tests, with apparently good small sample behavior. If the size of the literature on cointegration in the finite dimensional setting provides any indication, there remains enormous scope for further research on estimation, testing and forecasting with cointegrated functional time series.
	
		\appendix
		
		\section{Mathematical appendix}
		
			\subsection{Spectral properties of operator-valued functions}\label{aftapx}
			
			Consider a map \(A:U\to\mathcal{L}_\mathcal{H}\), where \(U\) is some open connected subset of \(\mathbb C\). We say that \(A\) is holomorphic on an open connected set \(D\subseteq U\) if, for each \(z_0\in D\), the limit
			\begin{equation}
			A^{(1)}(z_0)\coloneqq\lim_{z\to z_0}\frac{A(z)-A(z_0)}{z-z_0}
			\end{equation}
			exists in the norm of \(\mathcal L_\mathcal{H}\). It can be shown \citep[pp.\ 7--8]{GGK90} that holomorphicity on \(D\) in fact implies analyticity on \(D\), meaning that, for every \(z_0\in D\), we may represent \(A\) on \(D\) in terms of a power series
			\begin{equation}
			A(z)=\sum_{k=0}^\infty(z-z_0)^kA_k,\quad z\in D,
			\end{equation}
			where \(A_0,A_1,\ldots\) is a sequence in \(\mathcal L_\mathcal{H}\) not depending on \(z\). When the power series terminates, so that \(A(z)\) is polynomial in \(z\), then we say that \(A\) is an operator pencil.
			
			The set of points \(z\in U\) at which the operator \(A(z)\) is noninvertible is called the spectrum of \(A\), and denoted \(\sigma(A)\). The spectrum is always a closed set, and if \(A\) is analytic on \(U\), then \(A(z)^{-1}\) depends analytically on \(z\in U\setminus\sigma(A)\) \citep[p.\ 56]{M12}. A lot more can be said about \(\sigma(A)\) and the behavior of \(A(z)^{-1}\) if we assume that \(A(z)\) is a Fredholm operator for every \(z\in U\). In this case we have the following result, a proof of which may be found in \citet[pp.\ 203--204]{GGK90}. It is a crucial input to our main results.
			\begin{AFT}
				Let \(A:U\to\mathcal L_\mathcal{H}\) be an analytic Fredholm operator-valued function, and assume that \(A(z)\) is invertible for some \(z\in U\). Then \(\sigma(A)\) is at most countable and has no accumulation point in \(U\). Furthermore, for \(z_0\in\sigma(A)\) and \(z\in U\setminus\sigma(A)\) sufficiently close to \(z_0\), we have
				\begin{equation}\label{aftlaurent}
				A(z)^{-1}=\sum_{k=-m}^\infty(z-z_0)^kN_k,
				\end{equation}
				where \(m\in\mathbb N\) and \(N_{-m},N_{-m+1},\ldots\) is a sequence in \(\mathcal L_\mathcal{H}\) not depending on \(z\). The operator \(N_0\) is Fredholm of index zero and the operators \(N_{-m},\ldots,N_{-1}\) are of finite rank.
			\end{AFT}
			The analytic Fredholm theorem tells us that \(A(z)^{-1}\) is analytic except at a discrete set of points, which are poles. The technical term for this property of \(A(z)^{-1}\) is meromorphicity. For further reading on the spectral properties of operator-valued functions we suggest \citet{GGK90} and \citet{M12}.
		
		\subsection{Proofs}\label{appx}
		    
		    Here we provide proofs of all numbered results. The ordering of the proofs differs from the order in which the results were stated, so as to respect logical antecedence.
		    
			\begin{proof}[Proof of Theorem \ref{polethm}]
			It is obvious that \((3)\Rightarrow(4)\), because when we write \(\mathcal H=\alpha_1\oplus\Pi_1\beta_1^\perp\), we are asserting that \(\mathcal H=\alpha_1+\Pi_1\beta_1^\perp\) and that \(\alpha_1\cap\Pi_1\beta_1^\perp=\{0\}\). To establish the equivalence of the four conditions, it remains to show that \((4)\Rightarrow(1)\Rightarrow(2)\Rightarrow(3)\). Under Assumption \ref{assume1}(ii), \(\Phi(z)\) is analytic, Fredholm operator-valued, noninvertible at \(z=1\), and invertible at \(z=0\). The analytic Fredholm theorem therefore implies that \(\Phi(z)^{-1}\) is analytic on a punctured neighborhood \(D\) of \(z=1\) with a pole of order \(d\) at \(z=1\), and for \(z\in D\) admits a Laurent series around \(z=1\) as in \eqref{aftlaurentintro}. Further, \(\Phi(z)\) is analytic on \(D\cup\{1\}\) and thus for \(z\in D\) admits the Taylor series
			\begin{equation}\label{taylor}
			\Phi(z)=\sum_{k=0}^\infty\frac{1}{k!}\Phi^{(k)}(1)(z-1)^k.
			\end{equation}
			Combining \eqref{aftlaurentintro} and \eqref{taylor} we obtain, for \(z\in D\),
			\begin{align}
			\mathrm{I}&=\left(\sum_{k=-d}^\infty \Upsilon_k(z-1)^k\right)\left(\sum_{k=0}^\infty\frac{1}{k!}\Phi^{(k)}(1)(z-1)^k\right)\\
			&=\sum_{k=-d}^\infty\left(\sum_{j=0}^{d+k}\frac{1}{j!}\Upsilon_{k-j}\Phi^{(j)}(1)\right)(z-1)^k.\label{idexp}
			\end{align}
			
			Suppose that condition (1) is false, meaning that \(d>1\). Then the coefficients of \((z-1)^{-d}\) and \((z-1)^{-d+1}\) in the expansion of the identity in \eqref{idexp} must be zero. That is,
			\begin{align}
			\Upsilon_{-d}\Pi_0&=\zero,\label{N1zero}\\
			\Upsilon_{-d+1}\Pi_0+\Upsilon_{-d}\Pi_1&=\zero.\label{N2zero}
			\end{align}
			Equation (\ref{N1zero}) implies that \(\Upsilon_{-d}\alpha_1=\{0\}\), while equation (\ref{N2zero}) implies that \(\Upsilon_{-d}\Pi_1\beta_1^\perp=\{0\}\). If the condition (4) were valid, we could conclude that \(\Upsilon_{-d}=\zero\); however, this is impossible since \(\Upsilon_{-d}\) is the leading coefficient in the Laurent series \eqref{aftlaurentintro}, which is nonzero by construction. Thus if condition (4) is true then condition (1) must also be true: \((4)\Rightarrow(1)\).
			
			We next show that \((1)\Rightarrow(2)\). Suppose that (1) is true, meaning that \(d=1\). The coefficients of \((z-1)^{-1}\) and \((z-1)^0\) in the expansion of the identity in \eqref{idexp} must be equal to \(\zero\) and \(\mathrm{I}\) respectively. Since \(d=1\), this means that
			\begin{align}
			\Upsilon_{-1}\Pi_0&=\zero,\label{N1zero2}\\
			\Upsilon_{0}\Pi_0+\Upsilon_{-1}\Pi_1&=\mathrm{I}.\label{N2zero2}
			\end{align}
			It is apparent from \eqref{N2zero2} that \(\Upsilon_{-1}\Pi_1(x)=x\) for all \(x\in\beta_1^\perp\), and from \eqref{N1zero2} that \(\Upsilon_{-1}\mathrm{P}_{\alpha_1}=\zero\). We deduce that
			\begin{equation}\label{N2zero3}
			\Upsilon_{-1}\mathrm{P}_{\alpha_1^\perp}\Pi_1(x)=x\quad\text{for all }x\in\beta_1^\perp.
			\end{equation}
			This shows that \(\Lambda_1\) is left-invertible, hence injective. The spaces \(\alpha_1^\perp\) and \(\beta_1^\perp\) are of equal and finite dimension due to the fact that \(\Pi_0\) is Fredholm of index zero under Assumption \ref{assume1}(ii), so injectivity of \(\Lambda_1\) implies bijectivity. Thus we have shown that \((1)\Rightarrow(2)\).
			
			We next show that \((2)\Rightarrow(3)\). This amounts to showing that \((2)\Rightarrow(4)\) and that (2) implies
			\begin{equation}\label{trivialint}
			\alpha_1\cap \Pi_1\beta_1^\perp=\{0\}.
			\end{equation}
			Condition (2) implies that \(\mathrm{P}_{\alpha_1^\perp}\Pi_1\beta_1^\perp=\alpha_1^\perp\). We therefore have
			\begin{equation}
			\mathcal H=\alpha_1+\mathrm{P}_{\alpha_1^\perp}\Pi_1\beta_1^\perp=\alpha_1+(\mathrm{I}-\mathrm{P}_{\alpha_1})\Pi_1\beta_1^\perp.
			\end{equation}
			Since every element of \((\mathrm{I}-\mathrm{P}_{\alpha_1})\Pi_1\beta_1^\perp\) is the sum of an element of \(\alpha_1\) and an element of \(\Pi_1\beta_1^\perp\), this shows that every element of \(\mathcal H\) can be written as the sum of an element of \(\alpha_1\) and an element of \(\Pi_1\beta_1^\perp\), and so it is proved that \((2)\Rightarrow(4)\). To establish that condition (2) also implies \eqref{trivialint} we observe that any element \(y\in\alpha_1\cap\Pi_1\beta_1^\perp\) may be written as \(y=\Pi_1(x)\) for some \(x\in\beta_1^\perp\). Projecting both sides of this equality on \(\alpha_1^\perp\) gives \(0=\mathrm{P}_{\alpha_1^\perp}\Pi_1(x)\). The bijectivity of \(\Lambda_1\) asserted by condition (2) thus requires us to have \(x=0\), implying that \(y=0\). Thus \eqref{trivialint} is proved under condition (2), and we have shown that \((2)\Rightarrow(3)\).		
			
			It remains to verify that the residue \(\Upsilon_{-1}\) is as claimed when condition (1) is satisfied. We showed earlier in the proof that condition (1) implies \eqref{N2zero3}. We deduce from this that
			\begin{equation}\label{N2zero4}
			\Upsilon_{-1}\mathrm{P}_{\alpha_1^\perp}\Pi_1\Lambda_1^{-1}\mathrm{P}_{\alpha_1^\perp}=\Lambda_1^{-1}\mathrm{P}_{\alpha_1^\perp},
			\end{equation}
			or, more simply, \(\Upsilon_{-1}\mathrm{P}_{\alpha_1^\perp}=\Lambda_1^{-1}\mathrm{P}_{\alpha_1^\perp}\). We also showed earlier in the proof that condition (1) implies that \(\Upsilon_{-1}\mathrm{P}_{\alpha_1}=\zero\), or equivalently \(\Upsilon_{-1}\mathrm{P}_{\alpha_1^\perp}=\Upsilon_{-1}\) Thus \(\Upsilon_{-1}=\Lambda_1^{-1}\mathrm{P}_{\alpha_1^\perp}\), as claimed. 
			
		\end{proof}
		
		\begin{proof}[Proof of Corollary \ref{polethmp1}]
			Since \(\Pi_1=-\Phi_1\) and \(\Phi_1(x)=x\) for all \(x\in\beta_1^\perp\), we must have \(\Pi_1\beta_1^\perp=\beta_1^\perp\). The equivalence of conditions (1), (2) and (3) therefore follows from \hyperref[polethm]{Theorem \ref*{polethm}}.
			
			Obviously \((2)\Rightarrow(4)\). We will show that \((4)\Rightarrow(1)\) by showing that (4) implies bijectivity of \(\Lambda_1\), which was established in \hyperref[polethm]{Theorem \ref*{polethm}} to be necessary and sufficient for a simple pole. The operator \(\Lambda_1\) reduces in the case \(p=1\) to the restriction of \(-\mathrm{P}_{\alpha_1^\perp}\) to \(\beta_1^\perp\), so its kernel is \(\alpha_1\cap\beta_1^\perp\). Since the domain and codomain of \(\Lambda_1\) are spaces of equal and finite dimension, condition (4) therefore implies that \(\Lambda_1\) is bijective. Thus \((4)\Rightarrow(1)\).
			
			It remains to show that, when \(p=1\) and we have a simple pole, \(-\Lambda_1^{-1}\mathrm{P}_{\alpha_1^\perp}\) corresponds to projection on \(\beta_1^\perp\) along \(\alpha_1\). Clearly \(-\Lambda_1^{-1}\mathrm{P}_{\alpha_1^\perp}\) has kernel \(\alpha_1\) and range \(\beta_1^\perp\), so it remains only to show idempotency. Since \(-\Lambda_1\) is equal to the restriction of \(\mathrm{P}_{\alpha_1^\perp}\) to \(\beta_1^\perp\) when \(p=1\), we have
			\begin{equation*}
			(-\Lambda_1^{-1}\mathrm{P}_{\alpha_1^\perp})(-\Lambda_1^{-1}\mathrm{P}_{\alpha_1^\perp})=\Lambda_1^{-1}(-\Lambda_1)\Lambda_1^{-1}\mathrm{P}_{\alpha_1^\perp}=-\Lambda_1^{-1}\mathrm{P}_{\alpha_1^\perp},
			\end{equation*}
			as claimed.
		\end{proof}
		
		\begin{proof}[Proof of Theorem \ref{mainthm}]
			Under Assumption \ref{assume1}(ii), \(\Phi(z)\) is analytic, Fredholm operator-valued, noninvertible at \(z=1\), and invertible elsewhere in the closed unit disk. The analytic Fredholm theorem therefore implies that \(\Phi(z)^{-1}\) is analytic on an open disk centered at zero with radius exceeding one, except at the point \(z=1\), where it has a pole, which we have assumed to be simple. Let \(\Psi(z) = (1-z)\Phi(z)^{-1} \), defined at \(z=1\) by analytic continuation.
			
			The fact that \( \Psi(z)\) and \(\tilde{\Psi}(z)\) are analytic on an open disk centered at zero with radius exceeding one implies that the coefficients of their Taylor series around zero, \(\Psi(z)=\sum_{k=0}^\infty\Psi_kz^k\) and \(\tilde{\Psi}(z)=\sum_{k=0}^\infty\tilde{\Psi}_kz^k\), decay exponentially in norm. Under Assumption \ref{assume1}(i), the two series \(\sum_{k=0}^\infty\Psi_k(\varepsilon_{t-k})\) and \(\sum_{k=0}^\infty\tilde{\Psi}_k(\varepsilon_{t-k})\) thus converge in \(L^2_{\mathcal H}\), the latter validly defining \(\nu_t\in L^2_{\mathcal H}\). Applying the equivalent linear filters induced by \((1-z)\Phi^{-1}(z)\) and \(\Psi(z)\) to either side of the equality \(X_t-\sum_{j=1}^p\Phi_j(X_{t-j})=\varepsilon_t\), we find that
			\begin{equation}\label{MAI1}
			\Delta X_t=\sum_{k=0}^\infty\Psi_k(\varepsilon_{t-k}),\quad t\geq1,
			\end{equation}
			a moving average representation for \(\Delta X_t\). Moreover, since \( \Psi(z)  = -\Upsilon_{-1} + (1-z)\tilde{\Psi}(z)\), we may rewrite \eqref{MAI1} as
			\begin{equation} \label{diffeq1}
			\Delta X_t = -\Upsilon_{-1}(\varepsilon_t) + \Delta\nu_t, \quad t\geq1.
			\end{equation}
			Clearly, the process given by 
			\begin{align}
			X_0^\ast=\nu_0,\quad X_t^* = -\Upsilon_{-1}\left(\sum_{s=1}^t \varepsilon_{s}\right) + \nu_t, \quad t\geq1,
			\end{align}
			is a solution to the difference equation \eqref{diffeq1}. It is completed by adding the solution to the homogeneous equation \(\Delta X_t=0\), which is any time invariant \(Z_0\in L^2_{\mathcal H}\). Therefore, we obtain the representation \eqref{deltax}.
			
			Since \(\Upsilon_{-1}\) is the residue of \(\Phi(z)^{-1}\) at \(z=1\), it is apparent from the formula \(\Upsilon_{-1}=\Lambda_1^{-1}\mathrm{P}_{\alpha_1^\perp}\) given in Theorem \ref{polethm} that the range of \(\Upsilon_{-1}\) is \(\beta_1^\perp\). Moreover, the range of \(\Upsilon_{-1}\) is of positive and finite dimension by the analytic Fredholm theorem. This establishes claim (1).
			
			It remains to establish claim (2). Supposing that \(Z_0\) belongs to \(\beta_1^\perp\), if we take the inner product of a nonzero \(x\in\beta_1\) with both sides of \eqref{deltax} we obtain
			\begin{equation*}
			\langle x,X_t\rangle=\langle x,\nu_t\rangle=\sum_{k=0}^\infty\hat{\Psi}^x_k(\varepsilon_{t-k}),\quad t\geq1,
			\end{equation*}
			where the \(\hat{\Psi}^x_k\)'s are given by
			\begin{equation*}
			\hat{\Psi}^x_k(y)=\langle x,\tilde{\Psi}_k(y)\rangle,\quad y\in\mathcal H,
			\end{equation*}
			and decay exponentially in norm in \(\mathcal{L}_{\mathcal H, \mathbb C}\). (To see why, it may be helpful to recall Remark \ref{iplprem}.) The sequence of inner products \((\langle x,X_t\rangle,t\geq 1)\) is thus a standard linear process in \(\mathbb C\) with innovations in \(\mathcal H\). The innovation covariance operator \(\Sigma\) is positive definite under Assumption \ref{assume1}(i), and so \((\langle x,X_t\rangle,t\geq 1)\) is I(0) if \(\sum_{k=0}^\infty\hat{\Psi}^x_k\neq\zero\). Noting that
			\begin{equation*}
			\sum_{k=0}^\infty\hat{\Psi}^x_k(y)=\left\langle x,\sum_{k=0}^\infty\tilde{\Psi}_k(y)\right\rangle,\quad y\in\mathcal H,
			\end{equation*}
			we see that \(\sum_{k=0}^\infty\hat{\Psi}^x_k=\zero\) if and only if \(x\) belongs to the orthogonal complement to the range of \(\sum_{k=0}^\infty\tilde{\Psi}_k\). But \(\sum_{k=0}^\infty\tilde{\Psi}_k\) is the analytic part of the Laurent series of \(\Phi(z)^{-1}\) around \(z=1\) evaluated at \(z=1\), which is \(\Upsilon_0\). We showed in the proof of Theorem \ref{polethm} that \eqref{N2zero2} holds in the presence of a simple pole. From this equality we deduce that the sum of the ranges of \(\Upsilon_0\) and \(\Upsilon_{-1}\) is \(\mathcal H\). It is therefore impossible for a nonzero \(x\in\beta_1\) to belong to the orthogonal complement to the range of \(\Upsilon_0\), and we conclude that \(\sum_{k=0}^\infty\hat{\Psi}^x_k\neq\zero\). Thus \((\langle x,X_t\rangle,t\geq 1)\) is I(0).
			
			If we instead take the inner product of some \(x\notin\beta_1\) with both sides of \eqref{deltax} we obtain
			\begin{equation*}
			\langle x,X_t\rangle=\langle x,Z_0\rangle-\left\langle x,\Upsilon_{-1}\left(\sum_{s=1}^t\varepsilon_s\right)\right\rangle+\langle x,\nu_t\rangle,\quad t\geq1.
			\end{equation*}
			Differencing yields
			\begin{align*}
			\langle x,\Delta X_t\rangle&=-\langle x,\Upsilon_{-1}(\varepsilon_t)\rangle+\langle x,\Delta\nu_t\rangle\\
			&=-\hat{\Upsilon}^x_{-1}(\varepsilon_t)+\sum_{k=0}^\infty\hat{\Psi}^x_k(\varepsilon_{t-k})-\sum_{k=0}^\infty\hat{\Psi}^x_{k+1}(\varepsilon_{t-k-1}),\quad t\geq2,
			\end{align*}
			where the operator \(\hat{\Upsilon}_{-1}^x\in\mathcal L_{\mathcal H,\mathbb C}\) is defined by
			\begin{equation*}
			\hat{\Upsilon}^x_{-1}(y)=\langle x,\Upsilon_{-1}(y)\rangle,\quad y\in\mathcal H.
			\end{equation*}
			Thus we see that \((\langle x,\Delta X_t\rangle,t\geq 2)\) is a standard linear process in \(\mathbb C\) with innovations in \(\mathcal H\), with coefficients summing to \(-\hat{\Upsilon}_{-1}^x\). Since \(x\) does not belong to \(\beta_1\), the orthogonal complement to the range of \(\Upsilon_{-1}\), we know that \(\hat{\Upsilon}^x_{-1}\neq\zero\). Thus \((\langle x,\Delta X_t\rangle,t\geq 2)\) is I(0) and \((\langle x,X_t\rangle,t\geq 1)\) is I(1).
		\end{proof}
		
		\begin{proof}[Proof of Theorem \ref{polethm2}]
			It is obvious that \((3)\Rightarrow(4)\), so to establish the equivalence of the four conditions, we will show that \((4)\Rightarrow(1)\Rightarrow(2)\Rightarrow(3)\).
			
			To show \((4) \Rightarrow (1)\), suppose that (1) is false; we will deduce that then (4) must also be false. Applying the analytic Fredholm theorem in the same way as in the proof of Theorem \ref{polethm}, we obtain identity \eqref{idexp} for \(z\) in a punctured neighborhood \(D\subset U\) of \(z=1\). Here, \(d\) is the order of the pole of \(\Phi(z)^{-1}\) at \(z=1\), and \(\Upsilon_{-d}\neq\zero\). A simple pole is ruled out by assumption, while a pole of second order is ruled out since we are maintaining that (1) is not satisfied, so we must have \(d>2\). The coefficients of \((z-1)^{-d}\), \((z-1)^{-d+1}\)  and \((z-1)^{-d+2}\) in the expansion of the identity in \eqref{idexp} must therefore be zero. That is,
			\begin{align} 
			\Upsilon_{-d}\Pi_0&=\zero, \label{i2eq1}\\
			\Upsilon_{-d+1} \Pi_0+\Upsilon_{-d}  \Pi_1&=\zero,\label{i2eq2} \\
			\Upsilon_{-d+2} \Pi_0 + \Upsilon_{-d+1} \Pi_1 + \Upsilon_{-d}\Pi_2+\Upsilon_{-d}\Pi_1\Pi_0^+\Pi_1&=\zero. \label{i2eq3}
			\end{align}
			From \eqref{i2eq1} we have \(\Upsilon_{-d}\alpha_1=0\), implying that \(\Upsilon_{-d}\alpha_2=\Upsilon_{-d}\Pi_1\beta_1^\perp\). And from \eqref{i2eq2} we have \(\Upsilon_{-d}\Pi_1\beta_1^\perp=0\). Therefore,
			\begin{align}\label{shortprf1}
			\Upsilon_{-d}\alpha_2=0.	
			\end{align}
			From \eqref{i2eq3} we have
			\begin{align*}
			\Upsilon_{-d}\Pi_2\beta_2^\perp=[\Upsilon_{-d+1}+\Upsilon_{-d}\Pi_1\Pi_0^+]\Pi_1\beta_2^\perp.
			\end{align*}
			By composing both sides of \eqref{i2eq2} with \(\Pi_0^+\), and noting that \(\Pi_0\Pi_0^+\) coincides with the identity on \(\alpha_1\), we see that \(\Upsilon_{-d+1}+\Upsilon_{-d}\Pi_1\Pi_0^+\) is zero on \(\alpha_1\). Since \(\Pi_1\beta_2^\perp\) is contained in \(\alpha_1\), we deduce that
			\begin{align}\label{shortprf2}
			\Upsilon_{-d}\Pi_2\beta_2^\perp=0.	
			\end{align}
			If (4) were true then \eqref{shortprf1} and \eqref{shortprf2} would together imply that \(\Upsilon_{-d}=\zero\). But this is impossible because \(d\) is the order of our pole at \(z=1\) and the associated Laurent coefficient must be nonzero. Thus \((4)\Rightarrow(1)\).
			
			Next we show that \((1) \Rightarrow (2)\).  The coefficients of \((z-1)^{-2}\), \((z-1)^{-1}\) and \((z-1)^0\) in the expansion of the identity in \eqref{idexp} must be equal to \(\zero\), \(\zero\) and \(\mathrm{I}\) respectively. Suppose that (1) is true. Then equations \eqref{i2eq1}, \eqref{i2eq2} and \eqref{shortprf1} above hold with \(d=2\), and in place of \eqref{i2eq3} we have
			\begin{align} 
			\Upsilon_{0} \Pi_0+\Upsilon_{-1} \Pi_1 +\Upsilon_{-2}\Pi_2+\Upsilon_{-2}\Pi_1\Pi_0^+\Pi_1&=\mathrm{I}. \label{i22eq3}
			\end{align}
			It is apparent from \eqref{i22eq3} that
			\begin{align*}
			\Upsilon_{-2}\Pi_2(x)=x-[\Upsilon_{-1}+\Upsilon_{-2}\Pi_1\Pi_0^+]\Pi_1(x)\quad\text{for all }x\in\beta_2^\perp.
			\end{align*}
			As observed earlier, it follows from \eqref{i2eq2} that \(\Upsilon_{-1}+\Upsilon_{-2}\Pi_1\Pi_0^+\) is zero on \(\alpha_1\). Since \(\Pi_1\beta_2^\perp\) is contained in \(\alpha_1\), we deduce that \(\Upsilon_{-2}\Pi_2\) coincides with the identity on \(\beta_2^\perp\). Moreover, it follows from \eqref{shortprf1} that \(\Upsilon_{-2}\Pi_2=\Upsilon_{-2}\mathrm{P}_{\alpha_2^\perp}\Pi_2\). We conclude that \(\Upsilon_{-2}\mathrm{P}_{\alpha_2^\perp}\Pi_2(x)=x\) for all \(x\in\beta_2^\perp\). This shows that \(\Lambda_2\) is left-invertible, hence injective. The spaces \(\alpha_2^\perp\) and \(\beta_2^\perp\) are of equal and finite dimension under Assumption \ref{assume1}(ii) (recall Remark \ref{kercokerrem}), so injectivity of \(\Lambda_2\) implies bijectivity. Thus we have shown that \((1)\Rightarrow(2)\).
				
			It can be shown that \((2)\Rightarrow(3)\) by arguing exactly as we did in the proof of the corresponding implication in Theorem \ref{polethm}, but replacing all subscripts of \(1\) with a subscript of \(2\).
			
			It remains to verify our formulas for the Laurent coefficients \(\Upsilon_{-2}\) and \(\Upsilon_{-1}\). Formula \eqref{resdef2} for \(\Upsilon_{-2}\) may be obtained by arguing exactly as we did in the proof of the residue formula in Theorem \ref{polethm}, but replacing all subscripts of \(1\) with a subscript of \(2\). Formula \eqref{resdef22} for \(\Upsilon_{-1}\) takes a bit more work. The coefficients of \((z-1)^{-2}\), \((z-1)^{-1}\), \((z-1)^{0}\) and \((z-1)^{1}\) in the expansion of the identity in \eqref{idexp} must be equal to \(\zero\), \(\zero\), \(\mathrm{I}\) and \(\zero\) respectively. When \(d=2\), this yields equations \eqref{i2eq1}, \eqref{i2eq2} and \eqref{i22eq3} above, as well as
			\begin{align}\label{i22eq4}
			\Upsilon_{1}\Pi_0+\Upsilon_0\Pi_1+\Upsilon_{-1}[\Pi_2+\Pi_1\Pi_0^+\Pi_1]+\Upsilon_{-2}[\Pi_3+\Pi_1\Pi_0^+\Pi_1\Pi_0^+\Pi_1]&=\zero.
			\end{align}
			Our Hilbert space \(\mathcal H\) satisfies \(\mathcal H=\alpha_1\oplus(\alpha_2\cap\alpha_1^\perp)\oplus\alpha_2^\perp\). Thus \(\Upsilon_{-1}=\Upsilon_{-1}\mathrm{P}_{\alpha_1}+\Upsilon_{-1}\mathrm{P}_{\alpha_2\cap\alpha_1^\perp}+\Upsilon_{-1}\mathrm{P}_{\alpha_2^\perp}\). We shall proceed by obtaining expressions for the three terms on the right-hand side of this equality separately. As observed earlier, it follows from \eqref{i2eq2} that \(\Upsilon_{-1}+\Upsilon_{-2}\Pi_1\Pi_0^+\) is zero on \(\alpha_1\). The Moore-Penrose inverse \(\Pi_0^+\) is zero on \(\alpha_1^\perp\), so we have
			\begin{align}\label{piece1}
			\Upsilon_{-1}\mathrm{P}_{\alpha_1}&=-\Upsilon_{-2}\Pi_1\Pi_0^+.
			\end{align}
			Next we observe that \(\Upsilon_{-1}\Lambda_1\) is equal to the restriction of \(\Upsilon_{-1}\Pi_1-\Upsilon_{-1}\mathrm{P}_{\alpha_1}\Pi_1\) to \(\beta_1^\perp\). Consequently, \(\Upsilon_{-1}\) coincides with \(\Upsilon_{-1}\Pi_1\Lambda_1^+-\Upsilon_{-1}\mathrm{P}_{\alpha_1}\Pi_1\Lambda_1^+\) on the range of \(\Lambda_1\), which is \(\alpha_2\cap\alpha_1^\perp\), and so we have
			\begin{align}
			\Upsilon_{-1}\mathrm{P}_{\alpha_2\cap\alpha_1^\perp}&=\Upsilon_{-1}\Pi_1\Lambda_1^+\mathrm{P}_{\alpha_2\cap\alpha_1^\perp}-\Upsilon_{-1}\mathrm{P}_{\alpha_1}\Pi_1\Lambda_1^+\mathrm{P}_{\alpha_2\cap\alpha_1^\perp}\notag\\
			&=\Upsilon_{-1}\Pi_1\Lambda_1^+\mathrm{P}_{\alpha_1^\perp}-\Upsilon_{-1}\mathrm{P}_{\alpha_1}\Pi_1\Lambda_1^+\mathrm{P}_{\alpha_1^\perp},\label{working}
			\end{align}
			with the second equality following from the fact that \(\Lambda_1^+(\alpha_2^\perp\cap\alpha_1^\perp)=\{0\}\). It is apparent from \eqref{i22eq3} that
			\begin{align}\label{working2}
			\Upsilon_{-1}\Pi_1\Lambda_1^+&=[\mathrm{I}-\Upsilon_{-2}\Pi_2-\Upsilon_{-2}\Pi_1\Pi_0^+\Pi_1]\Lambda_1^+.
			\end{align}
			If we substitute into \eqref{working} our expressions for \(\Upsilon_{-1}\Pi_1\Lambda_1^+\) and \(\Upsilon_{-1}\mathrm{P}_{\alpha_1}\) appearing in \eqref{working2} and \eqref{piece1} respectively, we obtain
			\begin{align}\label{piece2}
			\Upsilon_{-1}\mathrm{P}_{\alpha_2\cap\alpha_1^\perp}&=[\mathrm{I}-\Upsilon_{-2}\Pi_2]\Lambda_1^+\mathrm{P}_{\alpha_1^\perp}.
			\end{align}
			
			We still require a suitable expression for \(\Upsilon_{-1}\mathrm{P}_{\alpha_2^\perp}\). Our formula \eqref{resdef2} reveals that the range of \(\Upsilon_{-2}\) is contained in \(\beta_1^\perp\), implying that \(\Pi_0\Upsilon_{-2}=\zero\). Therefore, by composing both sides of \eqref{i22eq4} with \(\Upsilon_{-2}\), we obtain
			\begin{align}\label{i22eq42}
			\Upsilon_0\Pi_1\Upsilon_{-2}+\Upsilon_{-1}[\Pi_2+\Pi_1\Pi_0^+\Pi_1]\Upsilon_{-2}+\Upsilon_{-2}[\Pi_3+\Pi_1\Pi_0^+\Pi_1\Pi_0^+\Pi_1]\Upsilon_{-2}&=\zero.
			\end{align}
			Our formula \eqref{resdef2} also reveals that the range of \(\Pi_1\Upsilon_{-2}\) is contained in \(\alpha_1\), implying that \(\Pi_0\Pi_0^+\Pi_1\Upsilon_{-2}=\Pi_1\Upsilon_{-2}\). Therefore, by composing both sides of \eqref{i22eq3} with \(\Pi_0^+\Pi_1\Upsilon_{-2}\), we obtain
			\begin{align}\label{sub}
			\Upsilon_0\Pi_1\Upsilon_{-2}&=[\mathrm{I}-\Upsilon_{-1}\Pi_1-\Upsilon_{-2}\Pi_2-\Upsilon_{-2}\Pi_1\Pi_0^+\Pi_1]\Pi_0^+\Pi_1\Upsilon_{-2}.
			\end{align}
			We may substitute the right-hand side of \eqref{sub} for \(\Upsilon_0\Pi_1\Upsilon_{-2}\) in \eqref{i22eq42}. After some simplification this yields
			\begin{align}\label{i22eq43}
			\Upsilon_{-1}\Pi_2\Upsilon_{-2}+\Pi_0^+\Pi_1\Upsilon_{-2}+\Upsilon_{-2}[\Pi_3-\Pi_2\Pi_0^+\Pi_1]\Upsilon_{-2}&=\zero.
			\end{align}
			Formula \eqref{resdef2} implies that \(\mathrm{P}_{\alpha_2^\perp}\Pi_2\Upsilon_{-2}=\mathrm{P}_{\alpha_2^\perp}\). We may therefore rewrite \eqref{i22eq43} as
			\begin{align}\label{i22eq44}
			\Upsilon_{-1}\mathrm{P}_{\alpha_2^\perp}=-\Upsilon_{-1}\mathrm{P}_{\alpha_2}\Pi_2\Upsilon_{-2}-\Pi_0^+\Pi_1\Upsilon_{-2}+\Upsilon_{-2}[\Pi_2\Pi_0^+\Pi_1-\Pi_3]\Upsilon_{-2}.
			\end{align}
			We may substitute the sum of the right-hand sides of \eqref{piece1} and \eqref{piece2} for \(\Upsilon_{-1}\mathrm{P}_{\alpha_2}\) in \eqref{i22eq44}. After some simplification this yields
			\begin{align}
			\Upsilon_{-1}\mathrm{P}_{\alpha_2^\perp}&=-[\Lambda_1^+\mathrm{P}_{\alpha_1^\perp}\Pi_2+\Pi_0^+\Pi_1]\Upsilon_{-2}\notag\\
			&\quad+\Upsilon_{-2}[\Pi_1\Pi_0^+\Pi_2+\Pi_2\Pi_0^+\Pi_1+\Pi_2\Lambda_1^+\mathrm{P}_{\alpha_1^\perp}\Pi_2-\Pi_3]\Upsilon_{-2}.\label{piece3}
			\end{align}
			Finally, if we sum the right-hand sides of \eqref{piece1}, \eqref{piece2} and \eqref{piece3}, we obtain the desired formula \eqref{resdef22} for \(\Upsilon_{-1}\).
		\end{proof}
				
		\begin{proof}[Proof of Theorem \ref{mainthm2}]
			Under Assumption \ref{assume1}(ii) we may apply the analytic Fredholm theorem to deduce that \(\Phi(z)^{-1}\) is analytic on an open disk centered at zero with radius exceeding one, except at the point \(z=1\), where it has a pole, which we assume here to be of second order. Let \(\Psi(z)=(1-z)^2\Phi(z)^{-1}\), defined at \(z=1\) by analytic continuation.
			
			As in the proof of \hyperref[mainthm]{Theorem \ref*{mainthm}}, the coefficients in the Taylor series \(\sum_{k=0}^\infty\Psi_kz^k\) and \(\sum_{k=0}^\infty\tilde{\Psi}_kz^k\) decay exponentially in norm, while \(\sum_{k=0}^\infty\Psi_k(\varepsilon_{t-k})\) and \(\sum_{k=0}^\infty\tilde{\Psi}_k(\varepsilon_{t-k})\) converge in \(L^2_{\mathcal H}\), the latter validly defining \(\nu_t\in L^2_{\mathcal H}\). Applying the equivalent linear filters induced by \((1-z)^2\Phi^{-1}(z)\) and \(\Psi(z)\) to either side of the equality \(X_t-\sum_{j=1}^p\Phi_j(X_{t-j})=\varepsilon_t\), we find that
			\begin{equation}\label{MAI2}
			\Delta^2 X_t=\sum_{k=0}^\infty\Psi_k(\varepsilon_{t-k}),\quad t\geq2,
			\end{equation}
			a moving average representation for \(\Delta^2 X_t\). Moreover, since \( \Psi(z)  = \Upsilon_{-2}-(1-z)\Upsilon_{-1} + (1-z)^2\tilde{\Psi}(z)\), we may rewrite \eqref{MAI2} as
			\begin{equation} \label{diffeq2}
			\Delta^2 X_t = \Upsilon_{-2}(\varepsilon_t)-\Upsilon_{-1}(\Delta\varepsilon_t) + \Delta^2\nu_t, \quad t\geq2.
			\end{equation}
			Clearly, the process given by 
			\begin{align}
			X_0^\ast=\nu_0,\quad X_t^* = \Upsilon_{-2}\left(\sum_{s=1}^t\sum_{r=1}^s \varepsilon_{r}\right)-\Upsilon_{-1}\left(\sum_{s=1}^t\varepsilon_{s}\right) + \nu_t, \quad t\geq1,
			\end{align}
			is a solution to the difference equation \eqref{diffeq2}. It is completed by adding the solution to the homogeneous equation \(\Delta^2X_t=0\), which is \(Z_0+tZ_1\) for any time invariant \(Z_0,Z_1\in L^2_\mathcal{H}\). Therefore, we obtain the representation \eqref{deltax2}.
						
			It was established in Theorem \ref{polethm2} that \(\Upsilon_{-2}=\Lambda_{2}^{-1}\mathrm{P}_{\alpha_2^\perp}\), so the range of \(\Upsilon_{-2}\) is \(\beta_2^\perp\). Moreover, the range of \(\Upsilon_{-2}\) is of positive and finite dimension by the analytic Fredholm theorem. This establishes claim (1).
			
			It remains to establish claim (2). Supposing that \(Z_0\) and \(Z_1\) belong to \(\beta_2^\perp\), if we take the inner product of a nonzero \(x\in\beta_2\) with both sides of \eqref{deltax2} we obtain
			\begin{align}\label{i2rep}
			\langle x,X_t\rangle&=-\left\langle x,\Upsilon_{-1}\left(\sum_{s=1}^t\varepsilon_s\right)\right\rangle+\langle x,\nu_t\rangle,\quad t\geq1.
			\end{align}
			As in the proof of Theorem \ref{mainthm}, \((\langle x,\nu_t\rangle,t\geq1)\) is a standard linear process in \(\mathbb C\) with innovations in \(\mathcal H\), and is not I(0) if and only if \(x\) belongs to the orthogonal complement to the range of \(\Upsilon_0\). We showed in the proof of Theorem \ref{polethm2} that \eqref{i22eq3} holds in the presence of a pole of second order. From this equality we deduce that the sum of the ranges of \(\Upsilon_0\), \(\Upsilon_{-1}\) and \(\Upsilon_{-2}\) is \(\mathcal H\). It is therefore impossible for a nonzero \(x\in\beta_2\) to be orthogonal to both of the ranges of \(\Upsilon_0\) and \(\Upsilon_{-1}\). Thus either the first term on the right-hand side of \eqref{i2rep} is zero or I(1) and the second is I(0) (if \(x\) is not orthogonal to the range of \(\Upsilon_{0}\)), or the first term is I(1) and the second is a standard linear process (if \(x\) is not orthogonal to the range of \(\Upsilon_{-1}\)). In either case, \((\langle x,X_t\rangle,t\geq1)\) is I(0) or I(1).
			
			If we instead take the inner product of some \(x\notin\beta_2\) with both sides of \eqref{deltax2}, and twice difference, we obtain
			\begin{align}\label{i2rep2}
			\langle x,\Delta^2 X_t\rangle&=\langle x,\Upsilon_{-2}(\varepsilon_t)\rangle-\langle x,\Upsilon_{-1}(\Delta\varepsilon_t)\rangle+\langle x,\Delta^2\nu_t\rangle,\quad t\geq3.
			\end{align}
			Similar to the proof of Theorem \ref{mainthm}, we thereby find that \((\langle x,\Delta^2 X_t\rangle),t\geq3)\) is a standard linear process in \(\mathbb C\) with innovations in \(\mathcal H\), with coefficients summing to \(-\hat{\Upsilon}_{-2}^x\), where
			\begin{equation*}
			\hat{\Upsilon}^x_{-2}(y)=\langle x,\Upsilon_{-2}(y)\rangle,\quad y\in\mathcal H.
			\end{equation*}
			Since \(x\) does not belong to \(\beta_2\), we know that \(\hat{\Upsilon}^x_{-2}\neq\zero\). Thus \((\langle x,\Delta^2 X_t\rangle,t\geq 3)\) is I(0) and \((\langle x,X_t\rangle,t\geq 1)\) is I(2).
		\end{proof}


\begin{thebibliography}{99}                                                                  
		
		\bibitem[Al Sadoon(2018)]{A18}
		\textsc{Al Sadoon, M.} (2018). The linear systems approach to rational expectations models. \textit{Econometric Theory}, \textbf{34}, 628--658.
		
		\bibitem[Beare(2017)]{B17}
		\textsc{Beare, B.\ K.} (2017). The Chang-Kim-Park model of cointegrated density-valued time series cannot accomodate a stochastic trend. \textit{Econ Journal Watch}, \textbf{14}, 133--137.
		
		\bibitem[Beare and Seo(2017)]{BS17}
		\textsc{Beare, B.\ K.\ and Seo, W.\ -K.} (2017). Representation of I(1) autoregressive Hilbertian processes. ArXiv e-print, arXiv:1701.08149v1 [math.ST].
		
		\bibitem[Beare, Seo and Seo(2017)]{BSS17}
		\textsc{Beare, B.\ K., Seo, J.\ and Seo, W.\ -K.} (2017). Cointegrated linear processes in Hilbert space. \textit{Journal of Time Series Analysis},  \textbf{38}, 1010--1027.
		
		
		\bibitem[Ben-Israel and Greville(2003)]{BG03}
		\textsc{Ben-Israel, A.\ and Greville, T.\ N.\ E.} (2003). \textit{Generalized Inverses: Theory and Applications}, 2nd ed. Springer, New York.
		
		
		\bibitem[Bosq(2000)]{B00}
		\textsc{Bosq, D.} (2000). \textit{Linear Processes in Function Spaces}. Springer, New York.
		
		\bibitem[Cerovecki and H\"{o}rmann(2017)]{CH17}
		\textsc{Cerovecki, C.\ and H\"{o}rmann, S.} (2017). On the CLT for discrete Fourier transforms of functional time series. \textit{Journal of Multivariate Analysis}, \textbf{154}, 282--295.
		
		\bibitem[Chan and Wei(1987)]{CW87}
		\textsc{Chan, N.\ H.\ and Wei, C.\ Z.} (1987). Asymptotic inference for nearly nonstationary AR(1) processes. \textit{Annals of Statistics}, \textbf{15}, 1050--1063.
		
		\bibitem[Chan and Wei(1988)]{CW88}
		\textsc{Chan, N.\ H.\ and Wei, C.\ Z.} (1988). Limiting distributions of least squares estimates of unstable autoregressive processes. \textit{Annals of Statistics}, \textbf{16}, 367--401.
		
		\bibitem[Chang, Hu and Park(2016)]{CHP16}
		\textsc{Chang, Y., Hu, B.\ and Park, J.\ -Y.} (2016). On the error correction model for functional time series with unit roots. Mimeo, Indiana University.
		
		\bibitem[Chang, Kim and Park(2016)]{CKP16}
		\textsc{Chang, Y., Kim, C.\ -S.\ and Park, J.\ -Y.} (2016). Nonstationarity in time series of state densities. \textit{Journal of Econometrics}, \textbf{192}, 152--167.
		
		\bibitem[Cheng and Phillips(2009)]{CP09}
		\textsc{Cheng, X.\ and Phillips, P.\ C.\ B.} (2009). Semiparametric cointegrating rank selection. \textit{Econometrics Journal}, \textbf{12}, S83--S104.
		
		\bibitem[Cheng and Phillips(2012)]{CP12}
		\textsc{Cheng, X.\ and Phillips, P.\ C.\ B.} (2012). Cointegrating rank selection in models with time-varying variance. \textit{Journal of Econometrics}, \textbf{169}, 155--165.
		
		
		\bibitem[Davidson, Hendry, Srba and Yeo(1978)]{DHSY78}
		\textsc{Davidson, J.\ E.\ H., Hendry, D.\ F., Srba, F.\ and Yeo, S.} (1978). Econometric modelling of the aggregate time-series relationship between consumers' expenditure and income in the United Kingdom. \textit{Economic Journal}, \textbf{88}, 661--692.
		
		
		\bibitem[Engle and Granger(1987)]{EG87}
		\textsc{Engle, R.\ F.\ and Granger, C.\ W.\ J.} (1987). Co-integration and error correction: representation, estimation and testing. \textit{Econometrica}, \textbf{55}, 251--276.
		
		\bibitem[Engsted and Johansen(1999)]{EJ99}
		\textsc{Engsted, T.\ and Johansen, S.} (1999). Granger's representation theorem and multicointegration. In Engle, R.\ F.\ and White, H.\ (Eds.), \textit{Cointegration, Causality and Forecasting: Festschrift in Honour of Clive Granger}, pp.\ 200-211. Oxford University Press, Oxford.
		
		\bibitem[Faliva and Zoia(2002)]{FZ02}
		\textsc{Faliva, M.\ and Zoia, M.\ G.} (2002). Matrix polynomials and their inversion: the algebraic framework of unit-root econometrics representation theorems. \textit{Statistica}, \textbf{62}, 187--202.
		
		\bibitem[Faliva and Zoia(2009)]{FZ09}
		\textsc{Faliva, M.\ and Zoia, M.\ G.} (2009). \textit{Dynamic Model Analysis: Advanced Matrix Methods and Unit-root Econometrics Representation Theorems}, 2nd ed. Springer, Berlin.
		
		\bibitem[Faliva and Zoia(2011)]{FZ11}
		\textsc{Faliva, M.\ and Zoia, M.\ G.} (2011). An inversion formula for a matrix polynomial about a (unit) root. \textit{Linear and Multilinear Algebra}, \textbf{59}, 541--556.
		
		\bibitem[Franchi(2007)]{F07}
		\textsc{Franchi, M.} (2007). The integration order of vector autoregressive processes. \textit{Econometric Theory}, \textbf{23}, 546--553.
		
		\bibitem[Franchi and Paruolo(2011)]{FP11}
		\textsc{Franchi, M.\ and Paruolo, P.} (2011). Inversion of regular analytic matrix functions: local Smith form and subspace duality. \textit{Linear Algebra and its Applications}, \textbf{435}, 2896--2912.
		
		\bibitem[Franchi and Paruolo(2016)]{FP16}
		\textsc{Franchi, M.\ and Paruolo, P.} (2016). Inverting a matrix function around a singularity via local rank factorization. \textit{SIAM Journal of Matrix Analysis and Applications}, \textbf{37}, 774--797.
		
		\bibitem[Franchi and Paruolo(2018)]{FP18}
		\textsc{Franchi, M.\ and Paruolo, P.} (2018). Cointegration in functional autoregressive processes. ArXiv e-print, arXiv:1712.07522v2 [econ.EM].
		
		\bibitem[Franchi and Paruolo(2019)]{FP19}
		\textsc{Franchi, M.\ and Paruolo, P.} (2019). A general inversion theorem for cointegration. \textit{Econometric Reviews}, in press.
		
		\bibitem[Gohberg, Goldberg and Kaashoek(1990)]{GGK90}
		\textsc{Gohberg, I., Goldberg, S.\ and Kaashoek, M.\ A.} (1990). \textit{Classes of Linear Operators, Vol.\ 1}. Birkh\"{a}user, Basel.
		
		\bibitem[Granger(1981)]{G81}
		\textsc{Granger, C.\ W.\ J.} (1981). Some properties of times series data and their use in econometric model specification. \textit{Journal of Econometrics}, \textbf{16}, 121--130.
		
		\bibitem[Granger(1983)]{G83}
		\textsc{Granger, C.\ W.\ J.} (1983). Cointegrated variables and error-correcting models. Mimeo, University of California, San Diego.
		
		\bibitem[Granger(1986)]{G86}
		\textsc{Granger, C.\ W.\ J.} (1986). Developments in the study of cointegrated economic variables. \textit{Oxford Bulletin of Economics and Statistics}, \textbf{48}, 213--228.
		
		\bibitem[Granger and Lee(1989)]{GL89}
		\textsc{Granger, C.\ W.\ J.\ and Lee, T.\ H.} (1989). Investigation of production, sales and inventory relationships using multicointegration and non-symmetric error correction models. \textit{Journal of Applied Econometrics}, \textbf{4}, S145--S159.
		
		\bibitem[Granger and Lee(1990)]{GL90}
		\textsc{Granger, C.\ W.\ J.\ and Lee, T.\ H.} (1990). Multicointegration. In Rhodes, G.\ F.\ and Fomby, T.\ B. (Eds.), \textit{Advances in Econometrics, Vol.\ 8: Co-integration, Spurious Regressions, and Unit Roots}, pp.\ 71--84. JAI Press, Greenwich.
		
		
		\bibitem[H\"{o}rmann and Kokoszka(2012)]{HK12}
		\textsc{H\"{o}rmann, S.\ and Kokoszka, P.} (2012). Functional time series. Ch.\ 7 in Rao, T.\ S., Rao, S.\ S.\ and Rao, C.\ R.\ (Eds.), \textit{Handbook of Statistics, Vol.\ 30: Time Series Analysis---Methods and Applications}, pp.\ 157--186. North-Holland, Amsterdam.
		
		\bibitem[Howland(1971)]{H71}
		\textsc{Howland, J.\ S.} (1971). Simple poles of operator-valued functions. \textit{Journal of Mathematical Analysis and Applications}, \textbf{36}, 12--21.
		
		\bibitem[Hu and Park(2016)]{HP17}
		\textsc{Hu, B.\ and Park, J.\ -Y.} (2016). Econometric analysis of functional dynamics in the presence of persistence. Mimeo, Indiana University.
		
		\bibitem[Johansen(1988)]{J88}
		\textsc{Johansen, S.} (1988). The mathematical structure of error correction models. \textit{Contemporary Mathematics}, \textbf{80}, 359--386.
		
		\bibitem[Johansen(1991)]{J91}
		\textsc{Johansen, S.} (1991). Estimation and hypothesis testing of cointegration vectors in Gaussian vector autoregressive models. \textit{Econometrica}, \textbf{59}, 1551--1580.
		
		\bibitem[Johansen(1992)]{J92}
		\textsc{Johansen, S.} (1992). A representation of vector autoregressive processes integrated of order 2. \textit{Econometric Theory}, \textbf{8}, 188--202.
		
		\bibitem[Johansen(1996)]{J96}
		\textsc{Johansen, S.} (1996). \textit{Likelihood-Based Inference in Cointegrated Vector Autoregressive Models}. Oxford University Press, Oxford.
		
		\bibitem[Johansen(2009)]{J09}
		\textsc{Johansen, S.} (2009). Representation of cointegrated autoregressive processes with application to fractional processes. \textit{Econometric Reviews}, \textbf{28}, 121--145.		
		
		
		\bibitem[La Cour(1998)]{L98}
		\textsc{La Cour, L.} (1998). A parametric characterization of integrated vector autoregressive (VAR) processes. \textit{Econometric Theory}, \textbf{14}, 187--199.
		
		\bibitem[Kurtz and Protter(1991)]{KP91}
		\textsc{Kurtz, T.\ G.\ and Protter, P.} (1991). Weak limit theorems for stochastic integrals and stochastic differential equations. \textit{Annals of Probability}, \textbf{19}, 1035--1070.
		
		\bibitem[Markus(2012)]{M12}
		\textsc{Markus, A.\ S.} (2012). \textit{Introduction to the spectral theory of polynomial operator pencils}. American Mathematical Society, Providence.
		 
		\bibitem[Nielsen, Seo and Seong(2019)]{NSS19}
		\textsc{Nielsen, M.\ \O., Seo, W.\ -K.\ and Seong, D.} (2019). Variance ratio test for the number of stochastic trends in functional time series. Queen's Economics Department Working Paper No.\ 1420.
		 
		\bibitem[Park and Phillips(1988)]{PkP88}
		\textsc{Park, J.\ -Y.\ and Phillips, P.\ C.\ B.} (1988). Statistical inference in regressions with integrated processes: Part I. \textit{Econometric Theory}, \textbf{4}, 468--497.
		
		\bibitem[Park and Phillips(1989)]{PkP89}
		\textsc{Park, J.\ -Y.\ and Phillips, P.\ C.\ B.} (1989). Statistical inference in regressions with integrated processes: Part II. \textit{Econometric Theory}, \textbf{5}, 95--131.
		
		\bibitem[Phillips(1986)]{P86}
		\textsc{Phillips, P.\ C.\ B.} (1986). Understanding spurious regressions in econometrics. \textit{Journal of Econometrics}, \textbf{33}, 311--340.
		
		\bibitem[Phillips(1988)]{P88}
		\textsc{Phillips, P.\ C.\ B.} (1988). Regression theory for near-integrated time series. \textit{Econometrica}, \textbf{56}, 1021--1043.
		
		\bibitem[Phillips(1991)]{P91}
		\textsc{Phillips, P.\ C.\ B.} (1991). Optimal inference in cointegrated systems. \textit{Econometrica}, \textbf{59}, 283--306.
		
		\bibitem[Phillips and Durlauf(1986)]{PD86}
		\textsc{Phillips, P.\ C.\ B.\ and Durlauf, S.\ N.} (1986). Multiple time series regression with integrated processes. \textit{Review of Economic Studies}, \textbf{53}, 473--495.
		
		\bibitem[Phillips and Hansen(1990)]{PH90}
		\textsc{Phillips, P.\ C.\ B.\ and Hansen, B.\ E.} (1990). Statistical inference in instrumental variables regression with I(1) processes. \textit{Review of Economic Studies}, \textbf{57}, 99--125.
		
		\bibitem[Phillips and Kheifets(2019)]{PK19}
		\textsc{Phillips, P.\ C.\ B.\ and Kheifets, I.} (2019). On multicointegration. Mimeo, Yale University.
		
		\bibitem[Phillips and Park(1988)]{PPk88}
		\textsc{Phillips, P.\ C.\ B.\ and Park, J.\ -Y.} (1988). Asymptotic equivalence of ordinary least squares and generalized least squares in regressions with integrated regressors. \textit{Journal of the American Statistical Association}, \textbf{83}, 111--115.
		
		\bibitem[Phillips and Solo(1992)]{PS92}
		\textsc{Phillips, P.\ C.\ B.\ and Solo, V.} (1992). Asymptotics for linear processes. \textit{Annals of Statistics}, \textbf{20}, 971--1001.
		
		
		\bibitem[Sargan(1964)]{S64}
		\textsc{Sargan, J.\ D.} (1964). Wages and prices in the United Kingdom: A study in econometric methodology. In Hart, P.\ E., Mills, G.\ and Whitacker, J.\ K.\ (Eds.), \textit{Econometric Analysis for National Economic Planning}, pp.\ 23--63. Butterworths, London.
		
		
		\bibitem[Schumacher(1991)]{S91}
		\textsc{Schumacher, J.\ M.} (1991). System-theoretic trends in econometrics. In Antoulas, A.\ C.\ (Ed.), \textit{Mathematical System Theory: The Influence of R.\ E.\ Kalman}. Springer, Berlin.
		
		
		\bibitem[Seo(2018)]{S18}
		\textsc{Seo, W.\ -K.} (2018). Cointegration and representation of integrated autoregressive processes in function space. ArXiv e-print, arXiv:1712.08748v3 [math.FA].
		
		\bibitem[Seo and Beare(2019)]{SB19}
		\textsc{Seo, W.\ -K.\ and Beare, B.\ K.} (2019). Cointegrated linear processes in Bayes Hilbert space. \textit{Statistics and Probability Letters}, \textbf{147}, 90--95.
		
		\bibitem[Steinberg(1968)]{S68}
		\textsc{Steinberg, S.} (1968). Meromorphic families of compact operators. \textit{Archive for Rational Mechanics and Analysis}, \textbf{31}, 372--379.
		
		\bibitem[Yoo(1987)]{Y87}
		\textsc{Yoo, B.\ -S.} (1987). Co-integrated time series: Structure, forecasting and testing. Doctoral thesis, University of California, San Diego.
		
	\end{thebibliography}
\end{document}